\documentclass[%
 reprint,
 amsmath,amssymb,
 aps,
]{revtex4-2}

\usepackage{graphicx}% Include figure files
\usepackage{dcolumn}% Align table columns on decimal point
\usepackage{bm}% bold math
\usepackage{tikz, tipa}
\usepackage{pgfplots}
\usepackage{cancel}
\usepackage{subcaption}
\usepackage{float}

\def\T{{\mathbb{T}}}

\usetikzlibrary{arrows.meta}
\tikzset{>={Latex[width=2mm,length=2mm]}}

\usepackage[francais]{babel}

\usepackage{amsmath, amsfonts, amssymb,latexsym}
\usepackage{wasysym,graphicx,stmaryrd,tikz,bclogo}
\usepackage{pgfplots,xcolor,amsbsy}

\addto\captionsfrench{}

\usepackage{amsmath, amsfonts, amssymb,latexsym}

\input epsf.sty

\newcommand{\BEQ}{\begin{equation}}     % Gleichungen Anfang ..
\newcommand{\BEA}{\begin{eqnarray}}
\newcommand{\BD}{\begin{displaymath}}
\newcommand{\EEQ}{\end{equation}}       % .. und Ende
\newcommand{\EEA}{\end{eqnarray}}
\newcommand{\ED}{\end{displaymath}}
\newcommand{\del}{\delta}
\newcommand{\Del}{\Delta}
\newcommand{\eps}{\varepsilon}          % epsilon

\renewcommand{\P}{\mathbb{P}}

\def\proba{{\mathbb{P}}}

\def\esper{{\mathbb{E}}}
\def\T{{\mathbb{T}}}

\def\Prim {\mathrm{Prim}}

\newcommand{\Medskip}{\medskip\noindent}
\newcommand{\Bigskip}{\bigskip\noindent}

\catcode`\@=11
\def\numberbysection{\@addtoreset{equation}{section}
        \def\theequation{\thesection.\arabic{equation}}}
\numberbysection

\begin{document}

\preprint{APS/123-QED}

\title{An introduction to random rule-based chemical networks}% Force line breaks with \\
%\thanks{A footnote to the article title}%

\author{Jérémie Unterberger}
\affiliation{%
IECL (Institut Elie Cartan de Lorraine), Université de Lorraine, Nancy, France
}%
\email{jeremie.unterberger@univ-lorraine.fr}

\date{\today}% It is always \today, today,
             %  but any date may be explicitly specified

\begin{abstract}
Large chemical networks appear in various branches of chemistry and biology, in particular, cellular 
metabolism and prebiotic chemistry. Detailed simulations of such networks are difficult, 
or even impossible for lack of kinetic data. Various strategies have been developed to
produce synthetic random networks mimicking the large scale organizational properties of experimental chemical networks. These random networks are however mathematical artefacts, 
which fail to reflect the general reactivity structure of chemistry. 

We present here a new class of random models of prebiotic (uncatalyzed) chemistry, based on 
context-independent rules, which is coherent with the general compositional logic of 
metabolism. The general organization of the random networks fits within the small world paradigm.   We get a phase diagram of the models through an approximate mapping to a solvable
tree growth model. Our predictions go beyond a purely abstract connectivity analysis of the reaction graph by
studying the diversity of chemical mechanisms, and singling out  evolutionary patterns,
such as autocatalysis and multistationarity, paving the road to   
possible open-ended evolution.

\end{abstract}

%\keywords{Suggested keywords}%Use showkeys class option if keyword
                              %display desired
\maketitle

\tableofcontents

%\renewcommand{\thefootnote}{\fnsymbol{footnote}}

%%%%%%%%%%%%%%%%%%%%%%
%%%%%%%%%%%%%%%%%%%%%%%

\section{General motivation}

%%%%%%%%%%%%%%%%%
%%%%%%%%%%%%%%%%%%

Understanding the chemical nature and general statistical features of large chemical reaction networks is a 
prominent question in Origin of Life studies. Historically, the first attempts to deal with
the inherent complexity have come from the purely theoretical side:
the theories of collectively autocatalytic sets (CASs), and reflexively autocatalytic,
food-generated sets (RAFs), see \cite{Kauf}, \cite{HorSte}, \cite{SteHor}, \cite{Hor} have been used extensively to explore autocatalysis in relation to abiogenesis, following ideas originally due to Eigen and 
Schuster \cite{Eigen}. The later appearance of large databases of chemical reactions, such as KEGG \cite{KEGG},  the metabolic network of 
Escherichia Coli \cite{EcoCyc}, or the abiotic reaction database \cite{Adam}, and the 
simultaneous development of bioinformatic tools, have made possible an exploration of large
networks of experimental origin, see e.g.  \cite{And0}, \cite{And1}, \cite{And2},  or  \cite{Peng} which  discusses hierarchical network organization and possible evolution in terms of 
 seed-dependent autocatalytic systems (SDAS), whereby an initial set of simple molecules (foodset) serves
 as initial generation of a multi-generation metabolic network obtained by recursively applying rules coming from 
 various databases. 

\Medskip On the computational side, advanced DFT techniques have allowed the development of software like 
ChemKin \cite{ChemKin} and RMG \cite{RMG} in combustion chemistry, or more versatile
algorithms like
SCINE/Chemoton \cite{UnsGriRei},  YARP \cite{ZhaSav} and others, with applications
to parametrization of the prebiotic chemistry space \cite{Wang,ZhaGarSav}. Random sampling methods
\cite{IsmChaHab} have also been used for a rapid exploration of the chemical space.  In parallel, ML techniques are 
aiding both the inference of effective kinetic rate laws from experiment and the computational exploration of chemical reaction networks, see e.g.  \cite{Mar,WenSpo}. While DFT techniques allow rapid screening of reaction pathways and potential energy surfaces, the computational burden associated with high-precision functionals and accurate transition-state searches creates a fundamental trade-off between accuracy and throughput. Consequently, comprehensive exploration of very large chemical spaces remains beyond reach. 

\Medskip  This  
has suggested the generation of synthetic random networks mimicking the large scale organizational properties of experimental chemical networks as an explorational tool, see \cite{MulFlaSta} 
and refs. within.  Graphs are typically sparse, but do not 
look like regular graphs, nor like the Erd\H{o}s-R\'enyi random graph model.   Several authors have tried to
fit these (hyper)graphs to various models of random graphs \cite{Pan}, \cite{Gho}.   Other families of 
models have also been introduced, such as  the  "Small World"
networks  of Watts and Strogatz \cite{WatStr}, or the
Albert-Barabasi model of preferential attachment \cite{AlbBar}. Metabolic networks are generally
highly clustered, despite being sparse, and have a short characteristic path length due to a small number
of 'hubs', i.e. metabolites with high connectivity, thus fitting qualitatively the small world paradigm
\cite{Rav}; a quantitative agreement, however, is not observed.

%%%%%%%%%%%%%%%%%%%%%%%%%ù

\Medskip The abstract nature of the above random models makes them equally applicable to describe networks of very different nature 
(social networks, gene networks...)  What is crucially missing here, from our point of view, is
some  chemically interpretable observable. Interesting questions from the chemist's and biologist's point of view, in a context of Origin of Life studies,  may be: what are the main mechanisms 
producing long molecules with carbon skeletons ? Are they sufficiently diverse to produce
a large enough set of biologically relevant molecules, but still selective enough to avoid 
the unhindered generation of the full set of organic molecules ? Are there any signatures of prebiotic dynamics, such as autocatalytic cycles or multistability, pointing to a 
possible proto-Darwinian evolution ? Is open-ended evolution to arbitrarily
complex mechanisms possible ? 

\Medskip Elementary chemical knowledge suggests that mechanisms are based on context-independent
assembly of functional groups.    The M\O D package 
\cite{MOD} makes possible a systematic generation of pathways based on a fixed, finite grammar, i.e.  non-random list of rules, instead of individual reactions, and stands in striking contrast
to Erd\"os-R\'enyi type models, wherein each edge is accepted independently from the others.

\Medskip We introduce here instead   new classes of random  network models, called {\em random rule-based chemical networks}, RRBCN for short. In a RRBCN,  {\em a primitive rule, if accepted, 
generates an infinite number of rules by context derivation}, see Fig. \ref{Fig:1} below 
for an illustration. Contrary to the M\O D package, we provide random acceptation rules for mechanisms, 
which generate random, possibly infinite grammars. This procedure circumvents the   task of 
enumerating all known mechanisms -- which is error-prone, subject to interpretation and
necessarily non-exhaustive --, and makes it possible to raise general questions about what 
{\em any} type of chemistry could be like, as a function of interpretable parameters.     We argue  that such models are robust, in the sense that dressing them with the main colors and specifics of biochemistry adds a few layers of complexity, but does not alter the general strategy and results outlined here, as discussed in the conclusion.

\Medskip The rest of the section is devoted to an informal construction of the model, put against a 
general chemical context, and 
exposition of main results.

%%%%%%%%%%%%%%%%%

\subsection{From chemical reactions to abstract rules}

%%%%%%%%%%%%%%

Synthesis chemistry is concerned with finding a synthesis path with best possible yield for a given
compound. Expert knowledge collected by chemists over the years has allowed the discovery of a huge number
 of reaction paths  producing millions of molecules; some of this knowledge has
been collected into databases, most of them specific to a given branch (pharmacology, metabolism, etc.), see e.g. \cite{KEGG} or \cite{PubChem}.

Because these databases are oriented towards synthesis goals, it is however fair to say that they fail
to consider chemical networks in a systematic way, i.e. they do not provide a systematic, orderly set of 
chemical rules for reactivity. From common undergraduate knowledge, it is widely acknowledged that
reactivity rules follow general mechanisms, based on the formation and breaking of bonds between
moieties or functional groups (i.e., molecule fragments), and that these rules are widely independent on the {\em context}, i.e.
on the group substituents (neighboring atoms), though kinetics, therefore yields, may depend on these to a considerable extent.

\Medskip The point of view here is that of {\em network expansion}. The prototypal example is that of
the famous '52 Miller-Urey experiment, which demonstrated  the synthesis of a large diversity (thousands of
 compounds) of  organic compounds from inorganic constituents (in gaseous state, methane, ammonia, hydrogen and water vapor, supposed at the time to be representative of the archaic atmosphere of Earth) in an origin of life scenario; but other experimental sets have been considered and are
 under close scrutiny, see e.g. \cite{Huck}. With modern-day computers, and using a chemical expertise on
 reaction mechanisms, one may contemplate describing 
 the distribution of products inductively in terms of generations; see \cite{Peng} for an example. Generation 0 is the set of organic
 constituents, often called {\em foodset}. The computer may consider all theoretical products of 
 a priori plausible reaction mechanisms
 involving generation 0 compounds; those experimentally detected make up generation 1. Then, for 
 every $n\ge 2$,  generation $n$
 compounds consists of all experimentally detected products of compounds of generations $\le n-1$.  

\Medskip The key notion here is that of {\em mechanisms} or {\em primitive reactions}. 
 Chemical rules, as envisioned in bioinformatics softwares such as \cite{MOD}, are  {\em context-independent}, as opposed to individual database reactions, see e.g. \cite{KEGG}. This is perfectly standard for chemists, but let us provide an example for non-chemists 
 (Fig. \ref{Fig:1}). An aldol reaction is any reaction of the above form, 
 where $R_1,R_2,R_3,R_4$ are placeholders for arbitrary groups. The numbering of carbon 
 (C) atoms shows that this is an addition reaction, where two fragments are coalesced by making a bond.  The underlying  primitive reaction is obtained by considering only the reaction between the red parts (equivalently, by assuming trivial placeholders), 
 \BEA &&  {\mathrm{(aldolization\ mechanism\ /\ primitive\ reaction)}} \nonumber \\ 
 && \qquad
 O=C_3 +  C_1-C_2=O  \to O-C_3-C_1-C_2=O. \nonumber\\
 \EEA
  Aldolizations are the main rules0 driving the multiple
 steps of the 
 autocatalytic formose cycle \cite{But, Bre}, and participates in the Krebs cycle and
 other main metabolic cycles of life.

\begin{figure}  
\centering
\begin{tikzpicture}[scale = 0.7]

\begin{scope}[shift={(0,0)}]
%\draw(-3,0) node {$\bar{\rho}:$};
\draw(0,0) node {$R_4 - {\color{red} C_3} - R_3$};
\draw[red](-0.05,0.3)--(-0.05,0.6); \draw[red](0.05,0.3)--(0.05,0.6);
\draw[red](0,0.9) node {$O$};

\draw[red](-0.5,-0.8) rectangle(0.5,1.2);
\end{scope}

\begin{scope}[shift={(6,0)}]
\draw(0,0) node {$R_1 - {\color{red} C_1 - C_2} -  R_2$};
\draw[red](0.55,0.3)--(0.55,0.6); \draw[red](0.65,0.3)--(0.65,0.6);
\draw[red](0.6,0.9) node {$O$};

\draw[red](-1,-0.8) rectangle(1,1.2);
\draw[red,->](-0.6,-0.4) arc(-30:-150:3.2 and 1.5);

\end{scope}

\draw(3,0) node {$+$};

\begin{scope}[shift={(3,-4)}]
\draw(-5,0) node {$\to$};
\draw[red](-2,-0.75) rectangle(2,1.25);
\draw(0,0) node {$R_4 - {\color{red} C_3 - C_1 - C_2} - R_2$};
\draw(0,0.6)--(0,1.5); \draw(0,1.8) node {$R_1$};
\draw(-1.5,-0.6)--(-1.5,-1.5); \draw(-1.5,-1.8) node {$R_3$};
\draw[red](1.2-0.05,0.3)--(1.2-0.05,0.6); \draw[red](1.2+0.05,0.3)--(1.2+0.05,0.6);
\draw[red](1.2,0.9) node {$O$};
\draw[red](-1.5,0.3)--(-1.5,0.6); \draw[red](-1.5,0.9) node {$O$};
\end{scope}
\end{tikzpicture}
\caption{Aldolization mechanism}
\label{Fig:1}
\end{figure}

\Medskip Addition reactions  involve some electron rearrangement, which often manifests itself in bond displacements (here, the double bond $O=C_3$ turned into $O-C_3$) or 
in charges. We simplify chemistry by
\begin{enumerate}
\item considering only addition rules and their reverses, fragmentation rules;
\item considering all bond types alike;
\item disregarding  reactants or products with charges or cycles;
\item not discussing in full detail branching reactants or products (see below).
\end{enumerate}
Molecular graphs have thus been turned into abstract words.  Direct translation of reactivity rules  such as the aldolization mechanism in this simplified language is not possible. However, coarse-graining functional groups such as $O=C$ or $O-C$ into an extra 
'atom' denoted $D$ suggests the simplified  word addition (coalescence) rule 
\BEQ R_4{\color{red}D_2} + {\color{red} CD_1}R_2 \longrightarrow R_4{\color{red} D_2CD_1}R_2 \label{eq:RCDR}
 \EEQ
which allows two free contexts $R_4,R_2$ only. Eq. (\ref{eq:RCDR}) is a rule 
derived (by context addition) from the primitive rule
\BEQ 
{\color{red}D_2} + {\color{red} CD_1} \longrightarrow {\color{red} D_2CD_1}  \label{eq:CD}
\EEQ
We tend to use the word {\em rule} instead of {\em reaction} in our model, though the two are used interchangeably. Note, however, that the primitive rule
 (\ref{eq:CD})  allows $R_4{\color{red} D_2}R_3$ to react
with $R_1 {\color{red} CD_1}R_2$ for arbitrary, non trivial $R_3,R_1$, yielding a branching product. 

\Medskip In the following, we describe simplified models of reactivity satisfying 
these requirements, where molecules have become words with letters in an alphabet 
$\cal A$ called atom set, and bonds have become word concatenation, denoted by a dot
$(\cdot)$ or simple juxtaposition; words are identified by the word reversal symmetry, e.g.  $CD = DC$, a remnant of the graph symmetries involved in the SMILES representation of molecules \cite{Wei}. We use upper case letters for atoms, and lower case letters for words, typically $w,x$ or $f$, possibly with lower indices.

%%%%%%%%%%%%%%%%%

\subsection{Anabolic and catabolic rules, complexity indices}  \label{subsection:complexity}

%%%%%%%%%%%%%%%

Metabolism is often organized into anabolic (synthetic) and catabolic (degradation)
subsystems, the best-known examples of which being the Calvin cycle of photosynthesis vs. the Krebs cycle of respiration. The  distinction between anabolic and catabolic reactions may have
been essential in primordial metabolism, as has been put forward with much elegance and persuasive power by \cite{BraSmi}. Anabolism allows the construction of long molecules
with carbon skeleton, in general through simple addition rules (including e.g. aldolizations), wherein one of the 
reactants is a short molecule provided by the environment. Catabolism, on the other hand,
cuts long molecules into short pieces by repeated fragmentations. 

Here we specialize anabolism by assuming addition rules of the type 
\BEQ  {\mathbf{(anabolic\ rule)}} \qquad  f+x\to f\cdot x,
\EEQ
 where
$f\in  {\cal F}$ is in   a fixed,  finite set of molecules called {\bf external set}, which
would correspond to generation 0 in the above SDAS formalism.  Catabolic rules are of the general fragmentation form 
\BEQ {\mathbf{(catabolic\ rule)}}\qquad  x_1\cdot x_2 \to x_1 + x_2.
\EEQ
 Both anabolic rules
and catabolic rules are a priori reversible (and detailed balance should hold). To avoid interference between the two sets of rules, we assume that
$x,x_1,x_2\not\in {\cal F}$. Anabolic and catabolic rules together make up a set of 
rules called {\bf admissible rules}.

\Medskip The case when catabolic rules are irreversible corresponds to 
the dilute regime of open chemical reaction networks described in \cite{UntNgh}, where foodset species $f$ are assumed to be
abundant external species, and other species $x$ are low-concentration internal species, making reverse reactions $x_1+x_2\to x_1\cdot x_2$ of kinetic order 2 negligible. This linear regime is important, for it allows a direct interpretation of autocatalysis as instability
of the trivial 0 fixed point, with dynamical growth rate equal to the highest eigenvalue (Lyapunov eigenvalue) of the generator of the dynamics.

\Medskip {\em Complexity index.} The number of covalent bonds involved in an aldolization is 4 (the 4 red bonds on the right-hand side
of Fig. 1); this means that letting the two reactants approach each other and form a bond involves a 
global reorganization of the associated electron orbitals. We say that the  {\em complexity 
index} of such a reaction is 4, independently of the substituents $R_1,R_2,R_3,R_4$; it is  equal to the 
complexity of the underlying primitive reaction stripped of the $R_i$. In our simplified model, we call 

\begin{center} $|x| = $ (number of atoms in $x) - 1$ \end{center}

\noindent the {\em level} of $x$; it is also equal to the number of edges of the associated graph. Then  the {\em complexity index} of a
primitive anabolic rule $\rho : f+x\to f\cdot x$, resp. catabolic rule $x\to 
x_1\cdot x_2$, is 
\BEQ m_{\rho} = |f|+|x|+1, \qquad {\mathrm{resp.}} \qquad m_{\rho} = |x|.
\EEQ
  The complexity index of rules derived from those is identical.

%%%%%%%%%%%%%%%%%%%%%%%%%%%

\subsection{RRBCN Models: general construction}

%%%%%%%%%%%%%%%%%%%

The kinetics of simple reactions not involving enzymes can in principle be modeled by mass-action kinetics, namely, the velocity of a reaction $\rho: X_1+\ldots+X_n\to X'_1+\ldots+X'_{n'}$ 
is of the form $k_{\rho}$ times the product $\prod_{i=1}^n [X_i]$ of the concentrations of the reactants; the constant $k_{\rho}$ is called {\em kinetic rate}.  Within a given experiment,
the reaction set may usually be split according to a finite number of rankings $(N_{\rho})_{\rho}$, from the 
fastest ones to the slowest ones. Those reactions that operate too slowly, typically because
their kinetic rate is too low, are simply discarded as irrelevant. For simplicity, we discuss
here only the case of an {\em acceptation/rejection scheme}:   
 we choose 
a cut-off kinetic parameter $k_{min}$ and let 
\BEQ N_{\rho} = \begin{cases} 1, \qquad k_{\rho}>k_{min} \\ 0 \qquad {\mathrm{else}} \end{cases} \EEQ
Reaction rules with $N_{\rho}=1$ are {\em accepted}, those with $N_{\rho}=0$ are {\em rejected}.
Here $(N_{\rho})_{\rho}$ are defined on the basis of  independent, but not necessarily identically distributed, Bernoulli random variables $\omega_{\rho}$; namely, $N_{\rho}=1$ if and only if  $\omega_{\underline{\rho}}=1$ for some reaction $\underline{\rho}$ from which $\rho$ may
be derived by adding context  (including the case when $\omega_{\rho}=1$).    The simplest possible models depend on $\rho$ only through the  complexity indices $m_{\rho}$ 
defined in \S \ref{subsection:complexity}; namely, we let $\P[\omega_{\rho} = 1]=p_{\rho}$,  where $p_{\rho} =  p(m_{\rho})$ is a function of $m_{\rho}$. This leads to the following
formal definition of a network $\bar{\cal R}$ drawn from a {\bf RRBCN (random rule-based chemical network) model}, 
\BEA &&  \bar{\cal R} = \{\rho\ {\mathrm{admissible}}\ |\ N_{\rho} = 1\}, \nonumber\\
 && \qquad  N_{\rho} \overset{!}{=}   \max_{\underline{\rho}\subset \rho}  \omega_{\rho}, \qquad \omega_{\rho} \sim  {\mathrm{Ber}}(p_{\rho}),  \nonumber\\
&& (\underline{\rho}\subset\rho) \Leftrightarrow  \begin{cases}  \underline{\rho}:\ 
f+x\to f\cdot x, \\ 
\qquad \rho: \ f+ x\cdot r \to f\cdot x\cdot r  \ {\mathrm{(anabolic)}} \\
 \underline{\rho}:\ x_1\cdot x_2 \to x_1 + x_2,\\ 
\ \  \rho: \ (r_1\cdot x_1)\cdot (x_2 \cdot r_2)
\to r_1\cdot x_1 + x_2\cdot r_2 \\ \qquad \qquad  {\mathrm{(catabolic)}}  \end{cases}   \nonumber\\
 \label{eq:omega_rho} 
\EEA
in terms of Bernoulli random variables,  where the max ranges over all admissible reactions $\underline{\rho}$ from which 
$\rho$ may
be derived by adding (possibly trivial) contexts $r,r_1,r_2$. 

\Medskip Such models differ from standard independent edge random graphs by a {\em completion
procedure}.  Namely, $\cal R$
be the set of admissible reactions $\rho$ such that $\omega_{\rho}=1$. Then the reaction
network,
$\overline{\cal R}$, is obtained from $\cal R$ by adding all context derivations of 
reactions in $\cal R$.

\Medskip {\em Random network generation.} Fix a maximum level $n_{max}$. For each 
admissible reaction $\rho:f+x\to f\cdot x$ or $x_1\cdot x_2\to x_1 + x_2$ with 
$x,x_1,x_2\not\in {\cal F}$ and $|f\cdot x|,|x_1\cdot x_2|\le n_{max}$, we generate an independent
variable $\omega_{\rho}$ as in (\ref{eq:omega_rho}). Then, for each admissible 
reaction $\rho$, we let $N_{\rho} \overset{!}{=} \max_{\underline{\rho}\subset \rho} 
\omega_{\underline{\rho}}$ as in (\ref{eq:omega_rho}). The resulting random network is denoted $\bar{\cal R} = 
\bar{\cal R}(n_{max})$. 

\noindent Instead of this two-stage process (global sampling 
of all Bernoulli variables, followed by completion by derivation), we use for simulations a more intuitive and much faster algorithm by induction on the level, which we describe for non-branching molecules. We construct at the same time a set Prim$(n)$ of  {\bf primitive rules} of the form $f+x\to f\cdot x$ or $x_1\cdot x_2\to x_1 + x_2$ with $|f\cdot x|=n$, resp. $|x_1\cdot x_2|=n$.  Once the reaction sets $\bar{\cal R}(n-1)$ and Prim$(n')$, $n'=0,\ldots,n-1$ have been generated,     
$\bar{\cal R}(n)$ is obtained from $\bar{\cal R}(n-1)$ by adding level $n$ rules
of two types:
\begin{enumerate}
\item {\bf (inductive completion by derivation)} adding derived rules 
$f + x\cdot r \to f\cdot x\cdot r$, $(r_1\cdot x_1)\cdot (x_2\cdot r_2) \to 
r_1\cdot x_1 + x_2 \cdot r_2$ for all primitive reactions $f+x\to f\cdot x$, resp. 
$x_1\cdot x_2\to x_1 + x_2$ with $|f\cdot x|<n$, resp. $|x_1\cdot x_2|<n$;
\item {\bf (new primitive rules)} defining Prim $(n)$ as those level $n$ reactions
$\rho: f+x\to f\cdot x$, resp. $x_1\cdot x_2\to x_1 + x_2$, {\em not} obtained by completion as in 1., such that 
$|f\cdot x|=n$, resp. $|x_1\cdot x_2|=n$, for which $\omega_{\rho}=1$ (decide by tossing a coin with 
P(head) = $p_{\rho}$).  
\end{enumerate} 

\Medskip We still need to define a Bernoulli parameter $p_{\rho}$ for each admissible rule $\rho$.

%%%%%%%%%%%%%%%%%%%%%%%
 
\subsection{RRBCN Models I and II} \label{subsection:RRBCN-I-II}

%%%%%%%%%%%%%%%%%%%%%%%%%

We consider in this work the following functional forms with three parameters $0<p,q,z<1$:

(1) {\bf (constant parameter model, or Model I)} $p_{\rho} \equiv p$ (anabolic),
$p_{\rho}\equiv q$ (catabolic reaction), where $0<p,q<1$ are constants; 

(2) {\bf (constant affinity parameter model, or Model II)}  $p_{\rho} \equiv p z^{m_{\rho}}$ (anabolic case), 
$p_{\rho} \equiv q z^{m_{\rho}}$ (catabolic case), for some parameters $0<p,q,z<1$. Parameter $z$ is
called {\em affinity} because it is conjugate to the complexity.

\Bigskip  {\em Model I} can be simply considered as random acceptation/rejection model
for primitive reactions. A 
possible motivation for {\em Model II} (from which Model I can be gotten back by letting $z\to 1$) is the following. {\em Model kinetic constants $k_m$ of  rules $\rho$ in 
 with  complexity index $m$ as} 
\BEQ k_m\sim k_{ref} \sigma^m  U,  \label{eq:km} \EEQ
 where $k_{ref}<k_{min}$ is some  low kinetic value,  $0<\sigma<1$ is a cross-section-like parameter, and $U$ is log-normal, namely, $U \sim 10^{\lambda Z}$ ($\lambda>0$), $Z\sim {\cal N}(0,1)$.  Then {\em accept $\rho$ if 
 \BEQ k_m > k_{min}  \label{eq:kmin} \EEQ
 for a certain kinetic threshold}
$k_{min}$, and let 
\BEQ p_{\rho} = q_{\rho}= \proba[k_m>k_{min}] \label{eq:p=q}
\EEQ

\Medskip   The idea behind the parameter $\sigma$ is
that requiring an electron to move across $m$ atoms has an average activation energy price roughly
proportional to $m$. The resulting rules can be mapped into a constant affinity parameter model in
a certain regime of parameters; the parameter dictionary $(\frac{k_{min}}{k_{ref}}, \sigma, \lambda)
\to (p,z)$ is spelled out in \S I. A. of Suppl. Info. 

\Medskip It makes sense to postulate different values for $p_{\rho}$ and $q_{\rho}$; see
\S \ref{subsection:connectivity}.

%%%%%%%%%%%%%%%%%%%%%%%%%%%

\subsection{Main results}

%%%%%%%%%%%%%%%%%%%%%%%%%%%%

We present in this work results of two very different types. 

\Medskip The first category of results  are network-based explorations based on simulations (\S 
\ref{section:low-level}). The main emphasis here is on general connectivity features of the random networks 
(including sparsity, clustering, path lengths) and potential dynamical biosignatures 
 (autocatalysis, multistability). Focusing on low-complexity mechanisms and low-level (short)
 molecules, typically up to 6-7 atoms, allows full-scale simulations. RRBCN networks are "Small-World" networks
 like extant metabolic networks, and their structure is rich enough to allow in principle autocatalysis and
 multistability. We say "in principle", because the direction (direct or reverse) of reactions, in- and out-fluxes of external  molecules, and ultimately, kinetic rates, are not determined by the model, though they
 dictate the dynamical features that actually take place. 
 
\Medskip The second category of results are explorations of high-complexity mechanisms mostly based on theory (\S \ref{section:high-level}). The theoretical input is based on an approximate mapping to a solvable
random tree growth model, which is particularly relevant to explore the high-level part of the networks. 
Partial simulations exploring long molecules (up to 30-50 atoms) prove excellent agreement with our predictions, but theory allows an exploration of regimes
which are completely inaccessible by simulations. On the whole, we answer many questions raised in the
Introduction relative to the diversity of high-complexity mechanisms, including the possibility of open-ended
evolution, within the framework of our model. 

\Medskip Understanding proto-evolutionary dynamics would require to put together our two categories of results -- notably, understand the organization of reaction subnetworks generated by the high-complexity 
mechanisms unveiled in  \S \ref{section:high-level}. Some very partial results in that direction are 
presented in Suppl. Info., see introduction to Main Text,  \S \ref{section:high-level} for more.

\Medskip A distinctive feature of our models is the comparatively small amount of independent random variables
needed to define them -- one per primitive rule --, because of automatic context
completion.  We get two limiting cases for Model II:
\begin{enumerate}
\item Letting $p,z\to 0$ leaves only a small number of primitive rules. Then most likely
the random graph is made up of small islands, and the average node connectivity is low.

\item Letting $p,z\to 1$ also leaves only a small number of primitive rules. However,
the random graph has most likely a high average node connectivity, with a giant connected component containing most molecules. 
\end{enumerate}
What is puzzling here is that both limits carry little diversity : all primitive rules
have a low complexity index. The two cases differ by the fact that  $\bar{\cal R}$ 
contains few 
level $n$ rules (compared to the number of admissible level $n$ rules, for $n$ large)
when $p,z\to 0$, as opposed to most   level $n$ rules when $p,z\to 1$. Thus most 
level $n$ rules in the latter
case are deterministically derived from low-complexity primitive rules.
 Exploring intermediate regimes shows radical departure from independent
edge random graph models, such as the Erd\"os-R\'enyi model.   
A detailed exploration of available mechanisms can be performed by analytical tools,  based on auxiliary trees called {\em composition} for anabolism, or {\em fragmentation trees} for 
catabolism. Models I and II behave very 
differently, with Model I exhibiting only a finite (albeit potentially large) number of 
mechanisms, while  there appears an infinite number of  arbitrarily complex rules in some  intermediate range of the Model II parameter $z$, suggesting the possibility of 
  open-ended evolution.

%%%%%%%%%%%%%%%%%%%%%%%%%%%%%%

\section{Exploration of low-level structure: connectivity and dynamics}  \label{section:low-level}

%%%%%%%%%%%%%%%%%%%%%%%%%

In this section, we explore the model by simulations which focus on small-level cut-off random networks 
$\bar{\cal R}(n_{max})$ with $n_{max} = 6$ or $7$ at most. The idea is to explore the
shape and dynamical properties of a random network produced by a number of
low-complexity mechanisms.

\Medskip We have performed two types of network simulations: equilibrium network
simulations (\S \ref{subsection:connectivity}), focusing on the connectivity
structure of the reaction graph; and out-of-equilibrium network simulations (\S 
\ref{subsection:autocata}), focusing
on dynamical features.

%%%%%%%%%%%%%%%%%%%%%%%

\subsection{Equilibrium network simulations : network connectivity}\label{subsection:connectivity}

%%%%%%%%%%%%%%%%%%%%%%%

We assume here that admissible rules are reversible and that rates follow mass action and satisfy detailed balance. 

\noindent {\em Equilibrium.} Thermodynamic equilibrium for the anabolic reaction $f+ x \overset{k}{\underset{\bar{k}}{\rightleftarrows}} f\cdot x$ commands that
$k /\bar{k}  = e^{-\Del G^0/RT}$, where $\Del G^0$ is the Gibbs energy of the reaction, and 
$(RT)^{-1}$ the temperature-dependent Gibbs parameter. We assume here a flat energy landscape 
for simplicity $(\Del G^0\simeq 0)$, so that $k\simeq \bar{k}$. Let $C_0$, resp. $C_1$ be a typical equilibrium concentration of external  $f$, resp. internal species $x$. Given that $[f]\sim C_0$,
the effective constant of the direct reaction is however $\sim C_0 \bar{k}$ (see below). 
Similarly, given a reversible catabolic reaction $x_1\cdot x_2 \overset{k}{\underset{\bar{k}}{\rightleftarrows}} x_1 + x_2$, we may assume that the effective constant of the reverse reaction is $\sim C_1 k$.  

\Medskip {\em Determination of $p,q$ parameters.} The above effective constants allow a desymmetrization of the two rule types.   Namely, letting $\rho$ be a level $m$ direct anabolic reactions,  we define $p_{\rho}$ by substituting the effective
constant to the original one, namely,   $k_{ref}\to C_0 k_{ref}$ in 
(\ref{eq:km}); equivalently,  $\rho$ is accepted if $k_m > k_{min}/C_0$, with $k_m$ as in 
(\ref{eq:km}). 
Similarly,  reverse level $m$  catabolic rules are accepted if  $k_m > k_{min}/C_1$. Reverse anabolic or direct catabolic rules obey the original threshold relation (\ref{eq:kmin}).  See 
Suppl. Info. \S I.B.  The scaleless parameter $C_1/C_0$  is an important feature of the model;  e.g. $C_1$  large and
$C_0$ small will favor reaction networks that are strongly prejudiced towards fragmentation (by reverse anabolic or direct catabolic rules), thus inhibiting the formation of long molecules.

\Bigskip The  plots below have been obtained by averaging over samples of equilibrium simulations. We present results for 
the substrate graph $\cal S$. By definition, $\cal S$-nodes are molecules, and its edges
are obtained as follows,

\begin{figure}[t]
\centering
\begin{tabular}{c|c}
 rule in $\bar{\cal R}$ & edges of $\cal S$  \\
\hline   $f+x\to f\cdot x$ & $(x,f\cdot x)$ \\
\hline   $x_1+ x_2 \to x_1\cdot x_2$ & $(x_1,x_1\cdot x_2), (x_2,x_1\cdot x_2)$ \\
\hline $f\cdot x \to f+x$ & $(f\cdot x,x)$ \\
\hline $x_1\cdot x_2\to x_1 + x_2$ & $(x_1\cdot x_2,x_1),(x_1\cdot x_2,x_2)$
  \end{tabular}
\caption{${\cal S}$-graph.}
\label{Fig:S-graph}
\end{figure}

\begin{figure*}[t]
\centering

\begin{tikzpicture}
\node at (0,0) {\includegraphics[scale=0.4]{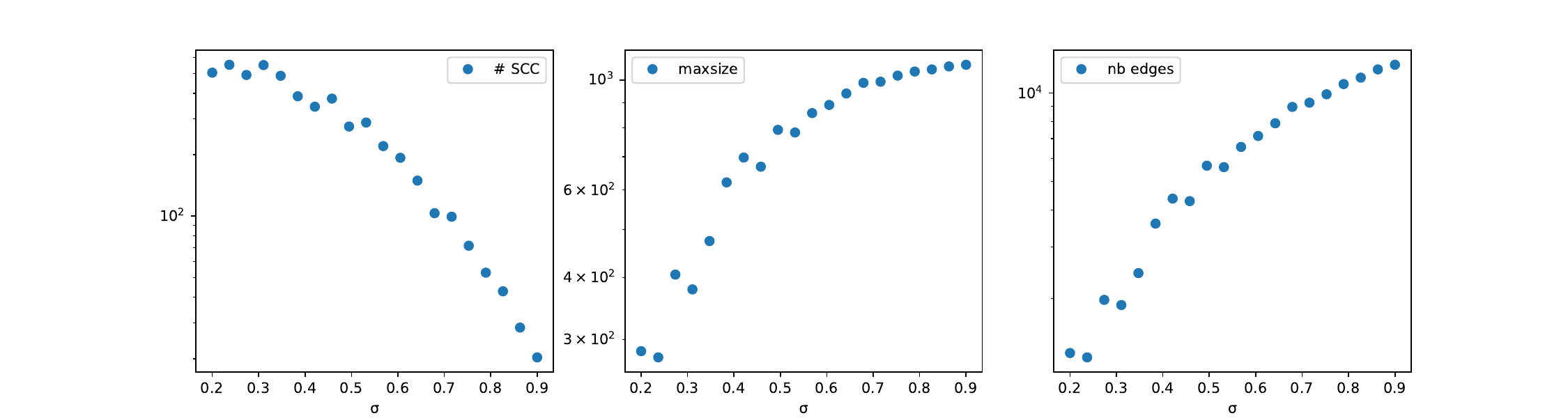}};
\node at (0,-5.5) {\includegraphics[scale=0.4]{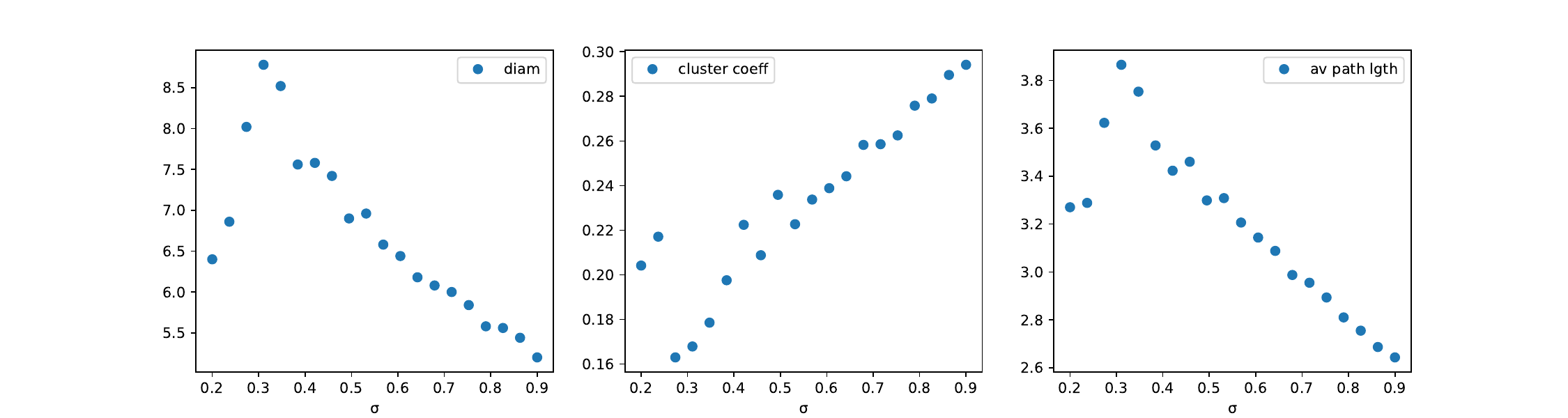}};

\node at (-4.14, 2.35) {\large\bfseries\textsf{(a)}};
\node at (-4.14+4.2, 2.35) {\large\bfseries\textsf{(b)}};
\node at (-4.14+8.4, 2.35) {\large\bfseries\textsf{(c)}};

\node at (-4.14, 2.85-6) {\large\bfseries\textsf{(d)}};
\node at (-4.14+4.2, 2.85-6) {\large\bfseries\textsf{(e)}};
\node at (-4.14+8.4, 2.85-6) {\large\bfseries\textsf{(f)}};

\end{tikzpicture}

\caption{Statistical features of $\bar{\cal R}(n_{max}=5)$. Above: (a) Number of SCCs. (b) Size of two largest SCCs. (c) Number of edges of 
largest SCC.  Below: considering only the largest SCC,  (d) Diameter. (e) Cluster coefficient. (f) Characteristic path length.}
\label{fig:ana_cata_size_and_cc_cpl_diam5}

\end{figure*}

\begin{figure*}[t]
\centering

\begin{tikzpicture}
\node at (0,0) {\includegraphics[scale=0.4]{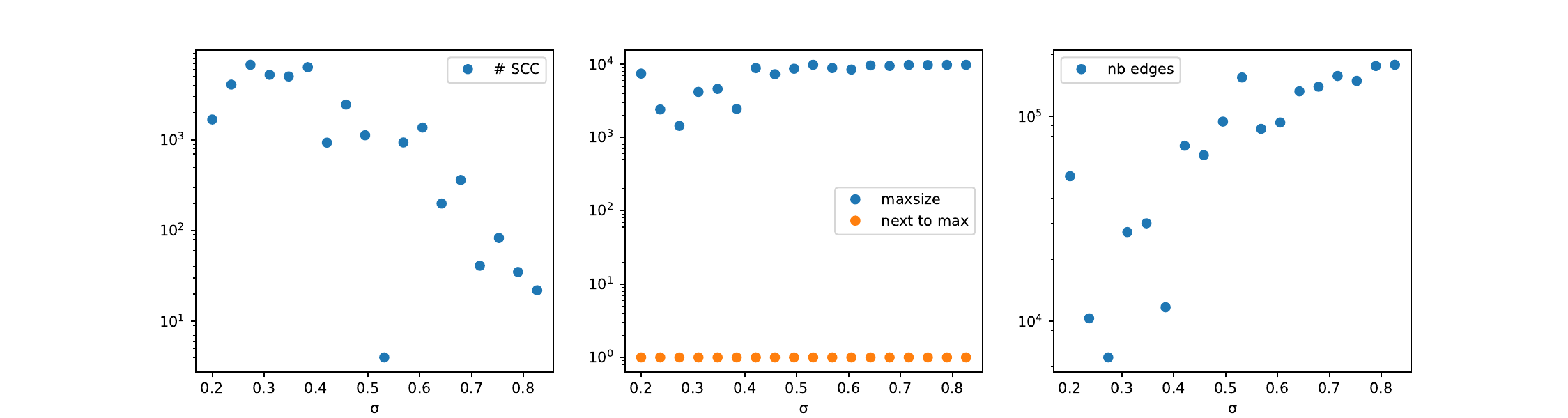}};
\node at (0,-5.5) {\includegraphics[scale=0.4]{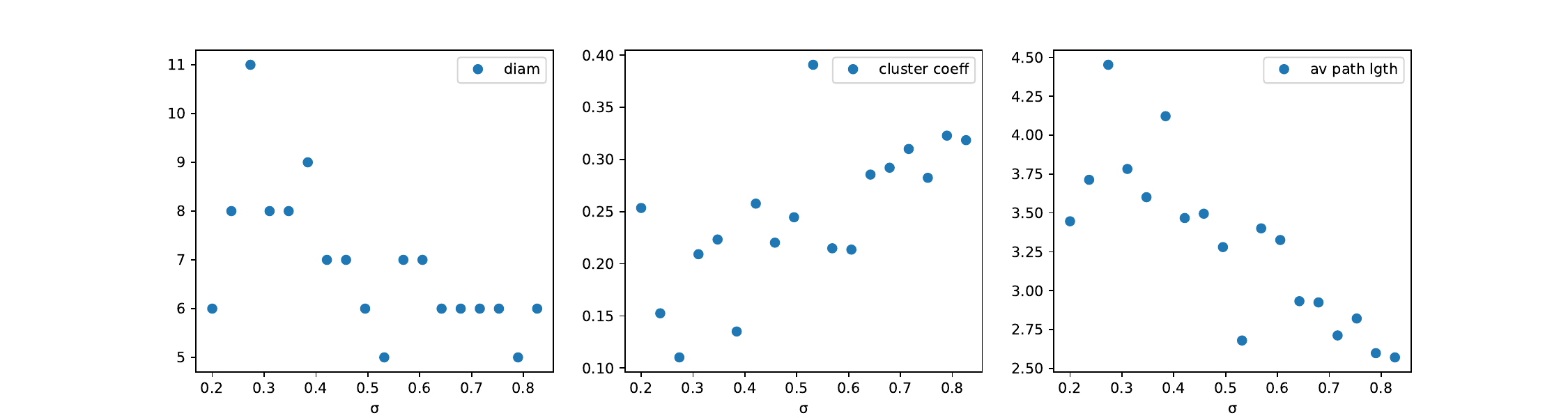}};

\node at (-4.14, 2.35) {\large\bfseries\textsf{(a)}};
\node at (-4.14+4.2, 2.35) {\large\bfseries\textsf{(b)}};
\node at (-4.14+8.4, 2.35) {\large\bfseries\textsf{(c)}};

\node at (-4.14, 2.85-6) {\large\bfseries\textsf{(d)}};
\node at (-4.14+4.2, 2.85-6) {\large\bfseries\textsf{(e)}};
\node at (-4.14+8.4, 2.85-6) {\large\bfseries\textsf{(f)}};

\end{tikzpicture}

\caption{Statistical features of $\bar{\cal R}(n_{max}=7)$. Above: (a) Number of SCCs. (b) Size of two largest SCCs. (c) Number of edges of 
largest SCC.  Below: considering only the largest SCC,  (d) Diameter. (b) Cluster coefficient. (c) Characteristic path length.}
\label{fig:ana_cata_size_and_cc_cpl_diam7}

\end{figure*}

\Medskip  We fix here
an alphabet ${\cal A}= \{A,B,C\}$, an external molecule set ${\cal F}=\{BC\}$ with only one element, a standard
deviation $\lambda=1$, $m_{ref}=\log_{10}(k_{min}/k_{ref}) = 2$, and a maximum level \BEQ n_{max}= \max\{|x|, x\in \bar{\cal R}\} = 5,
\EEQ
 see Fig. \ref{fig:ana_cata_size_and_cc_cpl_diam5} , and $n_{max}=7$, see Fig. \ref{fig:ana_cata_size_and_cc_cpl_diam7}
(up to 6 or 8 atoms), while letting the $\sigma$-parameter, see  eq. (1.3),  take 20 
values between 0.2 and 0.9.  Results for $n_{max}=7$ look very similar to those for $n_{max}=7$, except that they are less noisy, and  it took about 1 hour instead of a few minutes to generate the data.

%%%%%%%%%%%%%%%%%%%% hmax = 7 now

Fig. \ref{fig:ana_cata_size_and_cc_cpl_diam5} (a), (b), (c)  shows the number of SCCs, size $N(\sigma)$ (number of molecules) and number of edges
$E(\sigma)$  of the largest strongly connected component (SCC) of the substrate graph, average over a sample of 50 graphs.
Plot (b) also shows the size of the next-to-largest SCC, which is 1.  Unsurprisingly, both increase with $\sigma$, because the acceptance probability of any rule also does. It saturates
at  $1092 = 3+3^2 + \cdots + 3^6$ for $\sigma$ large enough, implying that all molecules of 
level $\le 5$ are in the same SCC. On the other hand, $E(\sigma)$, while large, does clearly not
scale like $N(\sigma)^2$, which means that the graph is very sparse. Comparison 
to results for $n_{max}=7$, see Fig. \ref{fig:ana_cata_size_and_cc_cpl_diam7}, suggests that 
$E(\sigma)\propto N(\sigma)$ instead, with a prefactor $\sim 10$; due to time limitations, we considered only one graph, which explains why data look noisy (all the more sore as $\sigma$ is small, since then only few primitive reactions are actually contributing).  Note that (contrary to the case of Erd\"os-R\'enyi graphs), no
phase transition is discernible. 

\Medskip Then Fig. \ref{fig:ana_cata_size_and_cc_cpl_diam5}  (d), (e), (f)  reports other quantities related to the largest SCC:
its diameter $D(\sigma)$ (maximum eccentricity), which is the largest graph distance (= shortest path length)
between two nodes; the clustering coefficient $0<C(\sigma)<1$; and its characteristic (average) path length
(average graph distance between two nodes), $L(\sigma)$. The fact that $C$ remains bounded 
below, and $L$ remains bounded from above by some small number, fits the small world paradigm. 
As $\sigma$ increases, the diameter decreases (despite the size increase), and so does the 
characteristic path length, because there appear shortcuts; at the same time, the clustering 
coefficient increases.

%%%%%%%%%%%%%%%%%%%%%%%%%%%%%%%%%%%%%

\subsection{Out-of-equilibrium simulations: dynamical features}
\label{subsection:autocata}

%%%%%%%%%%%%%%%%%%%%%%%%%%%%%%%%%%%%

Next, we kept the equilibrium assumption for anabolic rules, but dropped it for catabolic rules, turning them
into irreversible reactions, either $x_1\cdot x_2\to x_1 + x_2$ or $x_1 + x_2\to 
x_1\cdot x_2$. The  direction was chosen by favoring networks featuring a pair of
coupled autocatalytic cycles.    Namely, such networks exhibit multistationarity in some
regimes; see \S \ref{subsubsection:pair} for such an example of rule-based network. 
RRBCN networks are randomly generated rule-based networks, so that the network in 
\S \ref{subsubsection:pair} is an instance of RRBCN network. \S \ref{subsubsection:coupled} goes on prove that similar coupling structures between
autocatalytic cycles are very
common in RRBCN models, fostering the appearance of multiple stationary states.

%%%%%%%%%%%%%%%%%%%%

\subsubsection{A pair of coupled autocatalytic cycles} \label{subsubsection:pair}

%%%%%%%%%%%%%%%%%%%%

\Medskip The inspiration comes from a network with a pair of coupled autocatalytic cycles
which we now discuss. It is 
based on the following mechanisms (primitive rules) on a two-letter alphabet ${\cal A} = \{A,N\}$ with external molecule set ${\cal F} = \{A,NA\}$; for more readability, we use lower indices for repeated letters, 
e.g. $NA_4 N$ instead of  $NAAAAN$:   

(anabolism) \qquad  $A+ A_2  \underset{\bar{k}_0}{\overset{k_0}{\rightleftarrows}} A_3, 
\qquad NA + A_3 N  \underset{\bar{n}_0}{\overset{n_0}{\rightleftarrows}} NA_4N; $ 

(catabolism) \qquad $A_4 \overset{\bar{k}_1}{\to} A_2 + A_2, 
\ \ NA_4N \overset{\bar{n}_1}{\to} NA_2 + NA_2,\ \  A_2N + A_3 \overset{\kappa_1}{\to} A_2NA_3, \ \ A_2NA_3\overset{\kappa_2}{\to} A_2 + NA_3$.

  The network $\bar{\cal R}(n_{max})$ is the network obtained
from these rules
by completion, with cut-off at level $n_{max}$. It contains two autocatalytic cycles, 
denoted ${\cal C}_A$ and ${\cal C}_N$,
both of the same type as the core formose reaction network; the
autocatalytic cycle ${\cal C}_A$ for $n_{max}=3$, resp. ${\cal C}_N$ for $n_{max}=5$ takes the form $\{A_2 \overset{+A}{\rightleftarrows} A_3  
\overset{+A}{\rightleftarrows} A_4 \to 2A_2\}$, resp. $\{NA_2  \overset{+A}{\rightleftarrows} NA_3  
\overset{+AN}{\rightleftarrows} NA_4N \to 2NA_2\}$. The autocatalytic cycle ${\cal C}_A$ in $\bar{\cal R}(n_{max}=5)$ is extended to a cycle extending to species $A_5,A_6$; 
see e.g. \cite{Huck} for a visualization of a real-world formose network. Whatever the cut-off level, these cycles are known to be autocatalytic, so that compositions of the type $([A] = 1, [a] = \eps)$, resp. 
$([A]=[NA]=1, [n]=\eps)$ where  $a$, resp. $n$ is one of the species
of ${\cal C}_A$, resp. ${\cal C}_N$, and $\eps\ll 1$, are unstable, with the concentrations of all species in ${\cal C}_A$, resp. 
${\cal C}_N$ growing exponentially with time. We control the 
concentration of external molecules by a CSTR (continuously stirred chemical reactor) formalism, see e.g. \cite{Joshi} : 
$A$, resp. $NA$, is chemostatted via a high injection rate $i_A$, resp. $i_{NA}$, 
and degradated through a dilution rate $\del_A$, resp. $\del_{NA}$, yielding equilibrium
concentrations 
\BEQ [A]_{eq} = i_A/\del_A, \qquad [NA]_{eq} = i_{NA}/\del_{NA}. 
\EEQ
 This ensures that
concentrations of species within cycles saturate at a high level instead of growing
indefinitely. Since species in ${\cal C}_A$ and ${\cal C}_N$ are disjoint, this defines
two stable stationary states with compositions $[\cdot]_A$, resp. $[\cdot]_N$, with $([a]_A >0, [n]_A=0)$, resp. 
$([a]_N = 0, [n]_N>0)$ for any species $a$, resp. $n$ in ${\cal C}_A$, resp. 
${\cal C}_N$.

\Medskip The coupling of the two cycles  via the rules
$ A_2 N + A_3   \overset{\kappa_1}{\to} A_2NA_3 
\overset{\kappa_2}{\to} A_2 + NA_3 $ defines a 'cross-over' type sequence of reactions:
one $A$ atom is exchanged between the two reactants. (The two cycles can  be realized by taking the glyoxylate cycle for ${\cal C}_A$, and a nitrogenated version
for ${\cal C}_N$. Addition reactions are aldolizations, and the 'cross-over' is then a transamination. The actual chemical network
contains a few more intermediates).

\Medskip We present two sets of numerical simulations obtained using Julia's Catalyst and 
HomotopyContinuation libraries for chemical networks, which allow an easy input for
reactions, and a calculation of all stationary states:

\begin{enumerate}
\item one  for the restricted
${\cal C}_A(n_{max}=3)$ cycle coupled to ${\cal C}_N(n_{max}=5)$;

\item the other for
the unrestricted $\bar{\cal R}(n_{max}=5)$ network, coupling ${\cal C}_A(n_{max}=5)$
with ${\cal C}_N(n_{max}=5)$. 
\end{enumerate}

\begin{figure*}[htbp]
\centering

  \begin{subfigure}{0.48\textwidth}
        \centering
        \includegraphics[width=\linewidth]{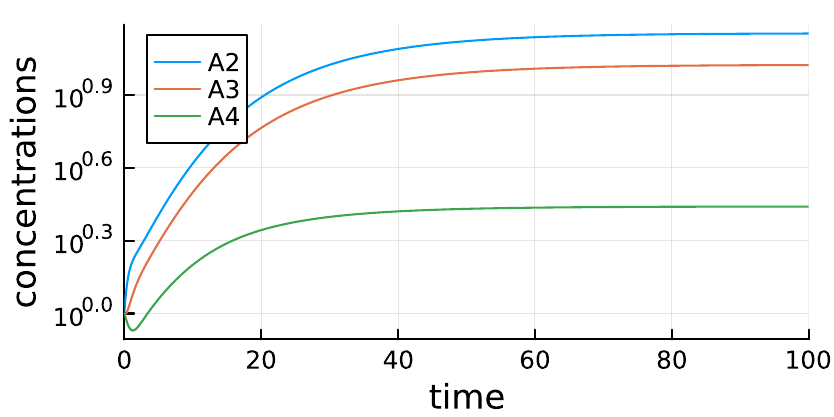}
        \caption{Total graph}
    \end{subfigure}
    \hfill
    \begin{subfigure}{0.48\textwidth}
        \centering
        \includegraphics[width=\linewidth]{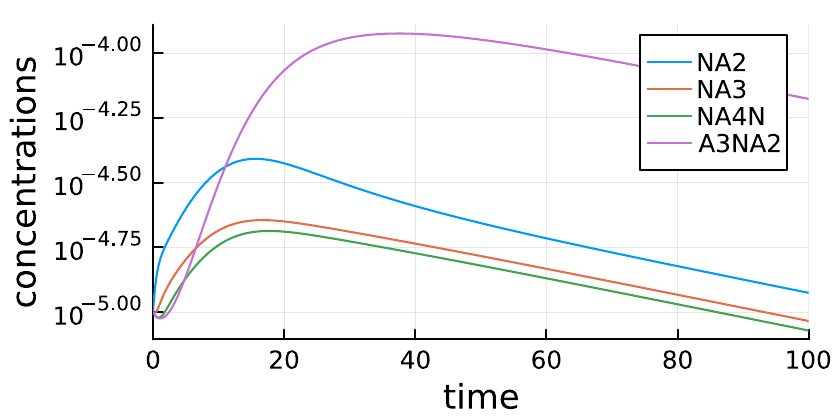}
        \caption{Connection graph}
    \end{subfigure}
\caption{2-cycle coupled network started from a neighborhood of the A-cycle}
\label{Fig:A}    
\end{figure*}

\begin{figure*}[htbp]
 \begin{subfigure}{0.48\textwidth}
        \centering
        \includegraphics[width=\linewidth]{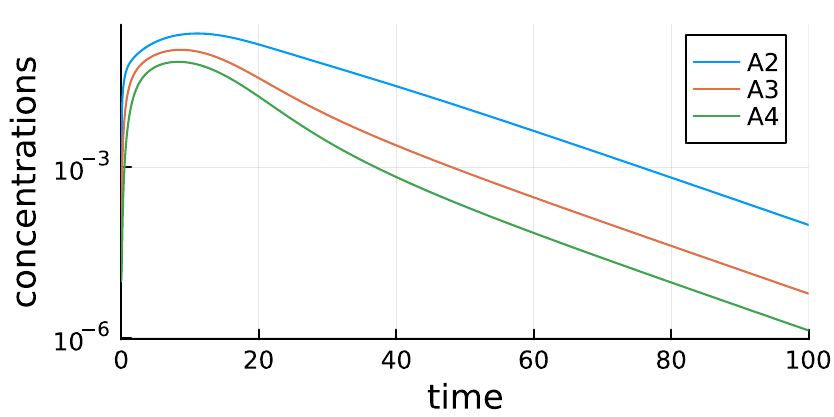}
        \caption{An autocatalytic SCC}
    \end{subfigure}
    \hfill
    \begin{subfigure}{0.48\textwidth}
        \centering
        \includegraphics[width=\linewidth]{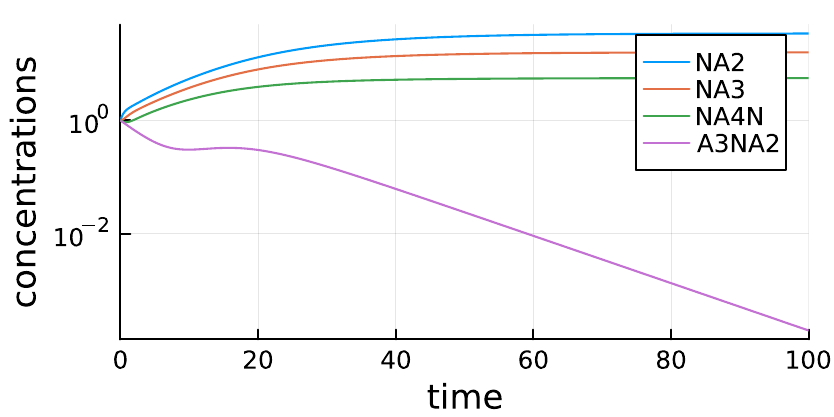}
        \caption{Another autocatalytic SCC}
    \end{subfigure}
\caption{2-cycle coupled network started from a neighborhood of the AN-cycle}
\label{Fig:NA}    
\end{figure*}

Parameter choices are specified in Suppl. Info. \S I.C.  Anabolic rules in case 1. have been
made irreversible, and rates of other reactions are in a close range.   In case 2., anabolic rules are reversible and detail-balanced with flat
thermodynamic balance $\Del G^0=0$, but rates are scattered over 8 ranges (from 
$\eps^{-3}$ to $\eps^5$, with $\eps = 1/4$). Dynamical simulations are presented in 
case 1. (see Fig. \ref{Fig:A}, \ref{Fig:NA}) which
start from the equilibrium values for $[A], [NA]$, and high values for ${\cal C}_A$
species/very low values for ${\cal C}_N$ species or vice-versa, corresponding to
possible invasion of ${\cal C}_A$ by ${\cal C}_N$ species or vice-versa. Simulations
show that invasions fail:  for instance, ${\cal C}_N$ species disappear in the first
case, while ${\cal C}_A$ species concentrations converge to their stationary values.
In short, single-cycle stationary states with compositions $[\cdot]_A$, $[\cdot]_N$ are stable.     

\begin{figure*}[htbp]
 %\begin{subfigure}{0.48\textwidth}
        \centering
        \includegraphics[scale=0.4]{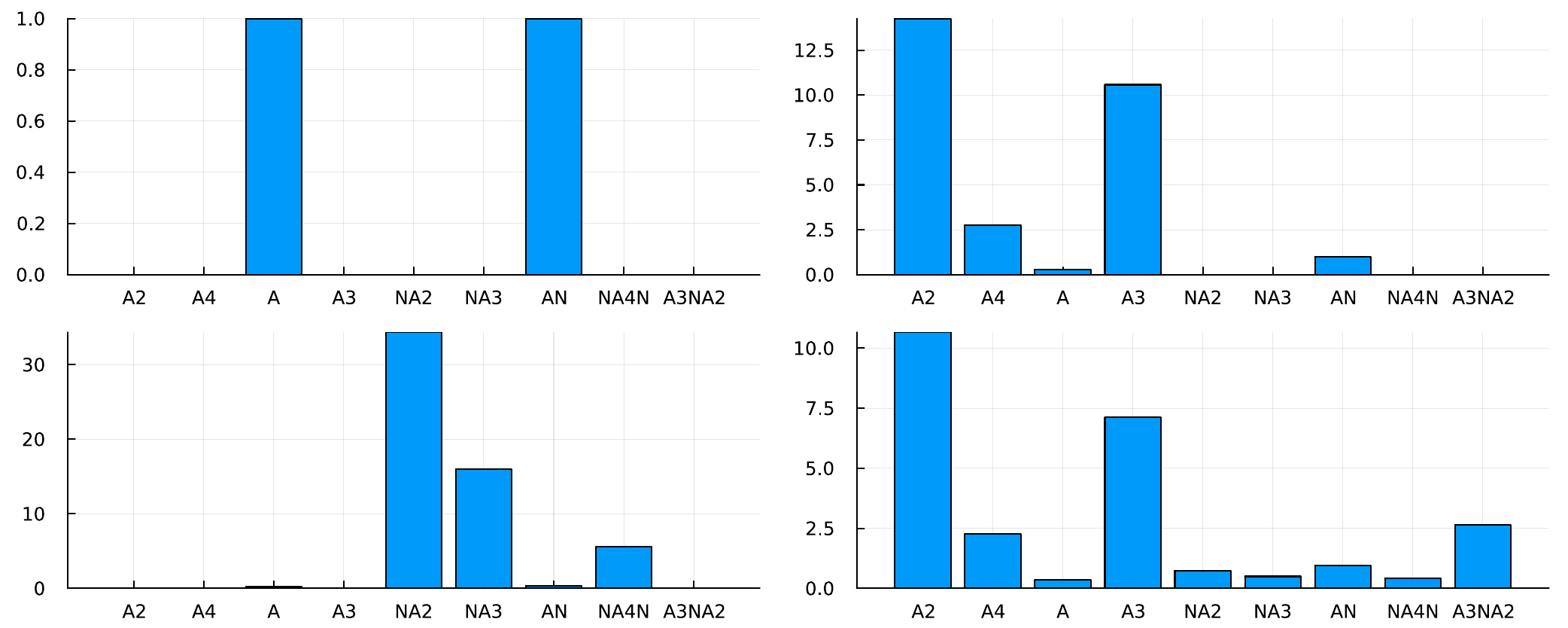}
        \caption{Stationary states for ${\cal C}_A(n_{max}=3)$ coupled to 
        ${\cal C}_N(n_{max}=5)$. From left to right and top to bottom:  $[\cdot]_{eq}, [\cdot]_A, [\cdot]_N, [\cdot]_{A,N}$. }
    %\end{subfigure}
    %\hfill
    %\begin{subfigure}{0.48\textwidth}
    %    \centering
    \label{Fig:2-cycle-stat35}  
\end{figure*}

\begin{figure*}[htbp]
        \includegraphics[scale=0.4]{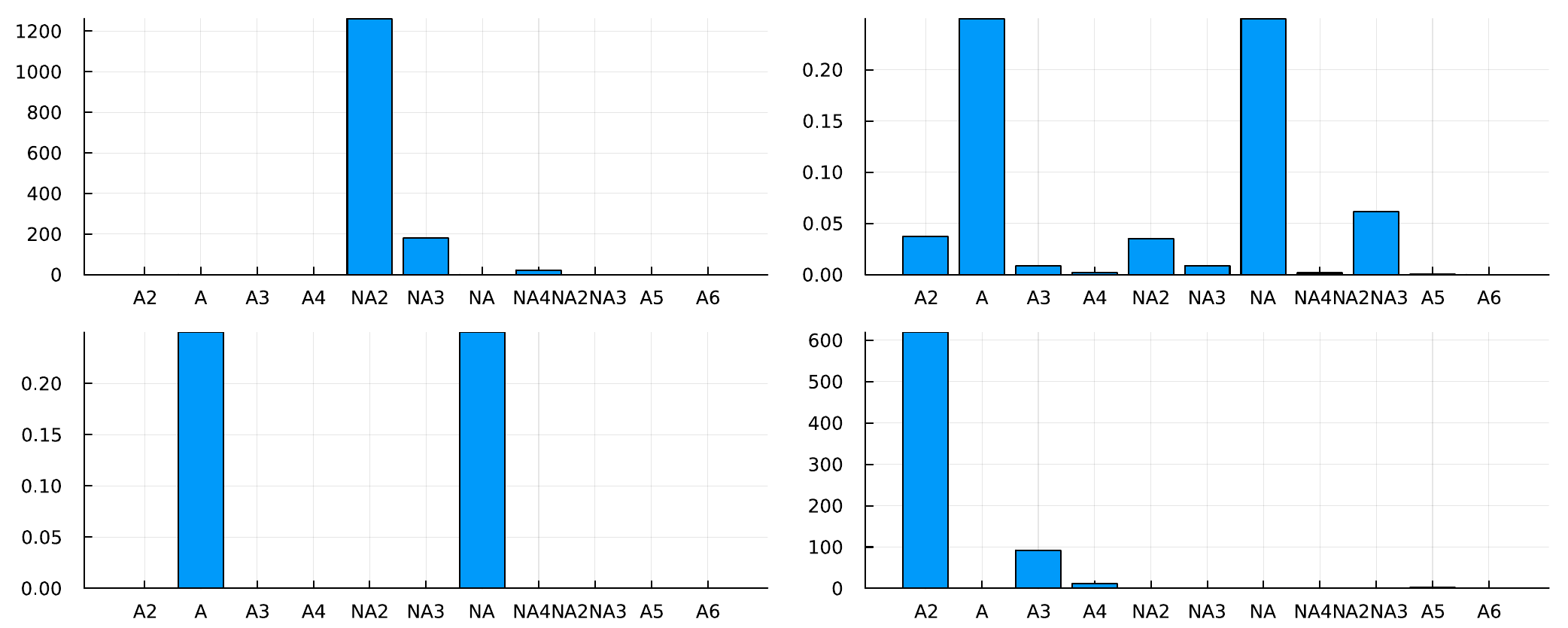}
        \caption{Stationary states  for ${\cal C}_A(n_{max}=5)$ coupled to 
        ${\cal A}_N(n_{max}=5)$.  From left to right and top to bottom:  $[\cdot]_{N}, [\cdot]_{A,N}, [\cdot]_{eq}, [\cdot]_{A}$. }
    %\end{subfigure}
%\caption{2-cycle stationary states}
\label{Fig:2-cycle-stat55}    
\end{figure*}

\Medskip In both cases 1. and 2., one gets 4 stationary states, see 
Fig. \ref{Fig:2-cycle-stat35} and \ref{Fig:2-cycle-stat55}: 

\begin{enumerate}
\item the two exclusive (single-cycle) stationary states with compositions 
$[\cdot]_A$, $[\cdot]_N$, which are stable;  plus two unstable states:
 
\item the trivial equilibrium  state with $[A]=[A]_{eq},\, [NA] = [NA]_{eq}$ and all
other concentrations vanishing;

\item and a cooperative state $[\cdot]_{A,N}$ (all A-cycle and N-cycle species present).  
\end{enumerate}

The cooperative state was proved to be unstable by checking that the Jacobian matrix
of the system at $[\cdot]_{A,N}$ has a positive eigenvalue; numerical simulations are 
a good illustration.

\vskip 0.5 cm


One can check that exclusive stationary states survive and remain stable when $\kappa_2$ is decreased down to $0$. Thus the general message is:  two autocatalytic cycles, coupled
by some reverse catabolic reaction $x_1+x_2\to x_1\cdot x_2$, with $x_1$ lying in 
one cycle, $x_2$ in the other, can produce stable exclusive stationary states. This generalizes in principle to multiple arrays of autocatalytic SCCs $\{C_i\}_{i=1,2,\cdots}$ coupled by catabolic rules $x_I + x_j \to x_i\cdot x_j$ defining an undirected graph
\textsf{ConnG} with nodes $\{1,2,\cdots\}$ and edges $\{(i,j)\ |\, x_i + x_j \to
x_i\cdot x_j\}$. Thus one expects to get a set of stable stationary states indexed
by independent sets  of \textsf{ConnG} (i.e. maximal subsets of \textsf{ConnG} satisfying
hard-core exclusion along edges)  as in  \cite{BunRiv}. We show in the
next paragraph how to generate such arrays in the RRBCN framework.

%%%%%%%%%%%%%%%%%%%%%%%

\subsubsection{Autocatalytic cycles and multistationarity in RRBCN models} \label{subsubsection:coupled}

%%%%%%%%%%%%%%%%%%%%%%%%%%%%%%%%%

Instead of picking a multistationary  network, we perform here an 
automatized search for  pairwise coupled autocatalytic cycles in RRBCN models. 
It turns out that this is a very common feature if the parameters of the models are
well-chosen.  The favorable parameter range  is difficult to assess a priori. We took here simply $p=q$ following (\ref{eq:km}, \ref{eq:kmin}), and chose to make all anabolic rules reversible.
It helped 
for numerics
that we picked $z=1$ (Model I) and $p=q=0.88$ high (see Suppl. Info. \S I.D.), so that the number
of rules is high even though the truncation level $n_{max}=5$ is low. We made 
simulations for two and three atom types; the results are similar, but graphs are
much larger for $|{\cal A}|=3$, so that they cannot be drawn properly with the Python networkx
package. Hence we kept to ${\cal A}=\{a,n\}$ as before, with lower-case letters this
time for convenience. The foodset here is ${\cal F} = \{a,n,na,an\}$ ($na$-additions are different
from $an$-additions because anabolic reactions are systematically written in the form 
$f+x\to f\cdot x$). Corresponding edges in the ${\cal S}$-graph are drawn red/blue/green/magenta if $f=a/n/na/an$. Edges coming from catabolic rules are black. 

All anabolic rules were chosen reversible.  We assumed for the network generation that catabolic 
reactions were irreversible (see Suppl. Info. \S I.D.). This amounts to neglecting reactions of the form 
$x_1 + x_2 \to x_1\cdot x_2$ which are of kinetic order 2. Assuming further that
foodset molecules $f$ are chemostated, we are in the "dilute" framework studied in 
\cite{UntNgh}, where the generator of the dynamics is linear, and a connectivity
analysis of the substrate graph $\cal S$ is enough to detect autocatalysis. Namely,
let ${\cal C}_i, i=1,2,\ldots$ be the strongly connected components (SCCs) of $\cal S$. 
Consider among the ${\cal C}_i$ those such that there exists a catabolic reaction 
$x_1\cdot x_2\to x_1 + x_2$ with all three species $x_1\cdot x_2, x_1, x_2$ in
${\cal C}_i$. Then (if the rates of rules exiting ${\cal C}_i$, i.e. with the reactant $x$
in ${\cal C}_i$ and at least one product outside, including dilution $x\to \emptyset$, are small enough), any trace of any 
species in ${\cal C}_i$ implies exponential growth of all species in ${\cal C}_i$, which
means that the network restricted to ${\cal C}_i$ is autocatalytic. We simply call
{\em autocatalytic SCCs} those ${\cal C}_i$ that satisfy the above topological property.

\Medskip Next, we consider the couplings between the autocatalytic SCCs. The conclusion
of \S \ref{subsubsection:pair} suggests that we look for sequences of  reverse catabolic reactions
coupling
${\cal C}_i$ and ${\cal C}_j$ $(i\not=j)$. We proceed as follows. Index all SCCs (autocatalytic or not) by
integer numbers; we first define a coarse-grained directed graph \textsf{G\_fork} on the set of SCC indices. Take any catabolic
rule $x\cdot y \rightleftarrows x + y$, and orient it in the reverse direction $x+y\to x\cdot y$. If $(x, x\cdot y)$  are in 
SCCs indexed by $(i,j)$, $i\not=j$,  draw an edge $i\to j$; similarly for the pair $(y,x\cdot y)$.    
Then we  build the undirected graph \textsf{ConnG} with: 
\begin{itemize}
\item[\textbullet] node set equal to the set of indices of autocatalytic SCCs; 
\item[\textbullet] edge set equal to the set of pairs   $(i,j)$ of nodes of \textsf{ConnG}
such that there exists a SCC $k$ (not necessarily autocatalytic) reachable both
from $i$ and from $j$ by a directed path in \textsf{G\_fork}.
\end{itemize} 
A particular case is when there exists a catabolic reaction $x+y\to x\cdot y$ with 
$x\in {\cal C}_i, y\in {\cal C}_j$, $x\cdot y\in {\cal C}_k$, $i\not=j\not=k$,  and $i,j \in $ \textsf{ConnG}: then
\textsf{G\_fork} contains the edges $(i,k),(j,k)$, so that $(i,j)$ are connected by an
edge in \textsf{ConnG}. 

\begin{figure*}[htbp]
\centering

  \begin{subfigure}{0.46\textwidth}
        \centering
        \includegraphics[width=\linewidth]{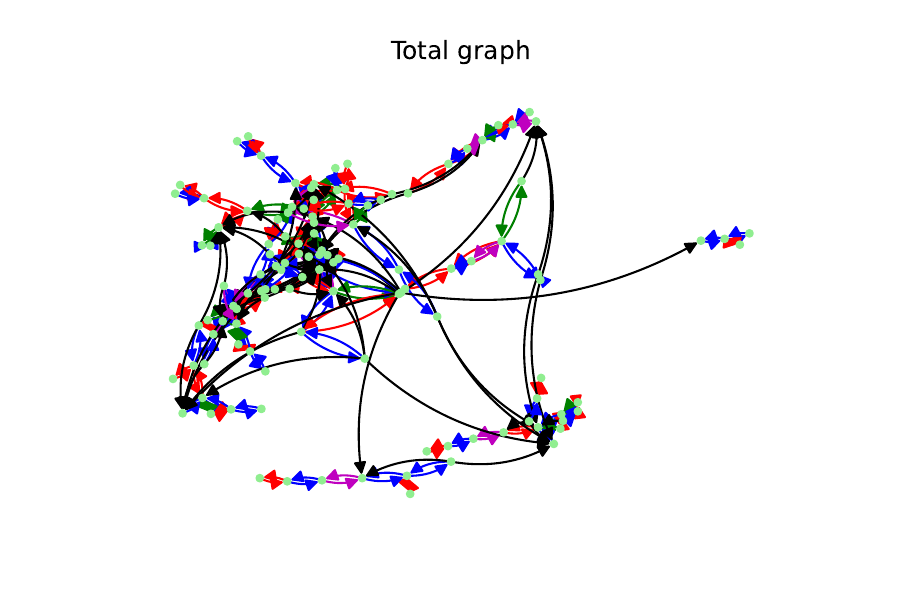}
        \caption{Total graph}
    \end{subfigure}
    \hfill
    \begin{subfigure}{0.46\textwidth}
        \centering
        \includegraphics[width=\linewidth]{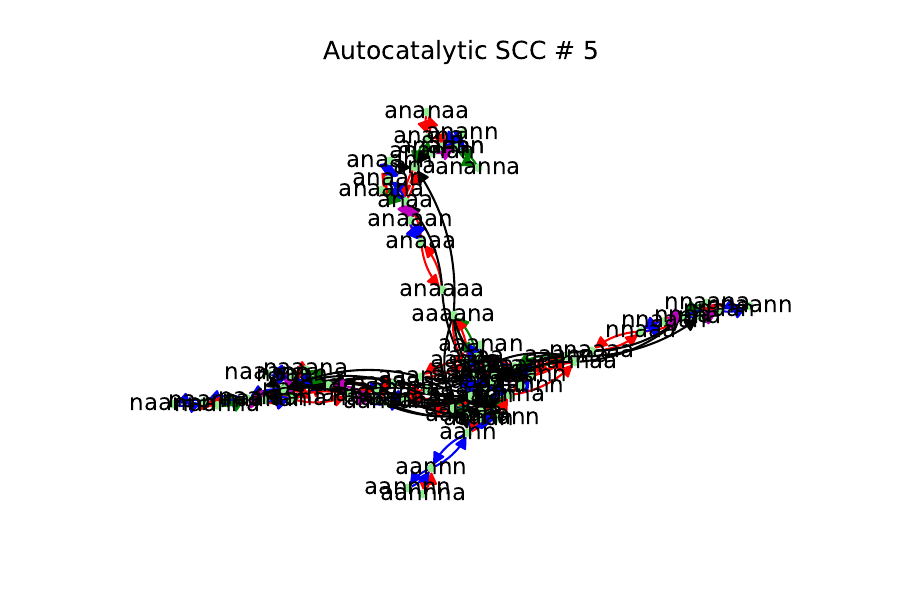}
        \caption{An autocatalytic SCC}
    \end{subfigure}
    
\vspace{0.5cm}

 \begin{subfigure}{0.48\textwidth}
        \centering
        \includegraphics[width=\linewidth]{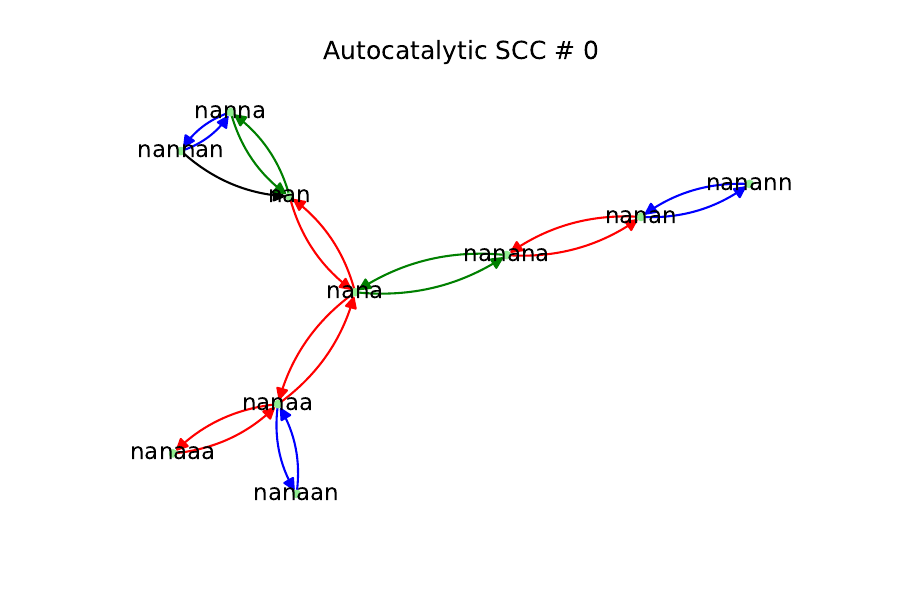}
        \caption{Another autocatalytic SCC}
    \end{subfigure}
    \hfill
    \begin{subfigure}{0.48\textwidth}
        \centering
        \includegraphics[width=\linewidth]{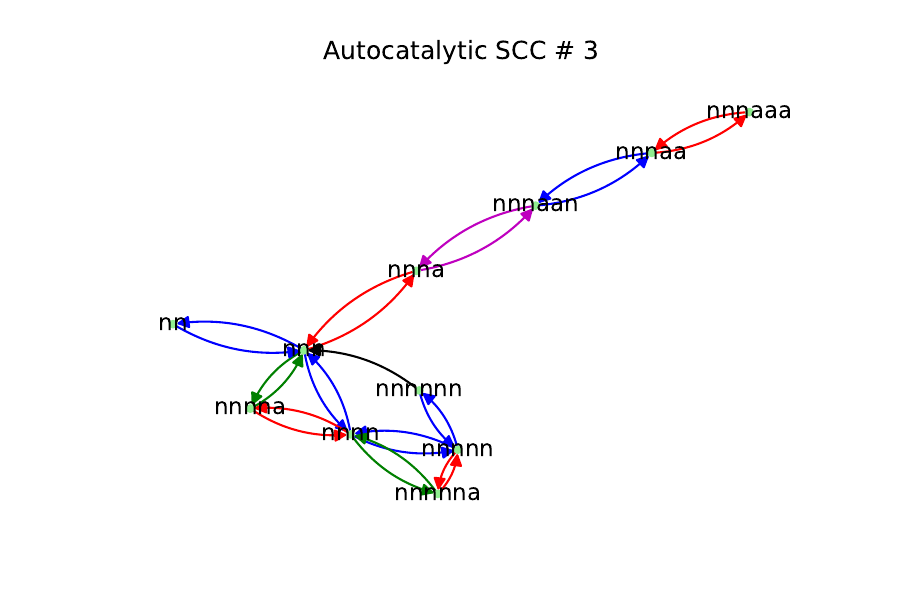}
        \caption{Yet another autocatalytic SCC}
    \end{subfigure}

\vspace{0.5cm}

\begin{subfigure}{0.48\textwidth}
        \centering
        \includegraphics[width=\linewidth]{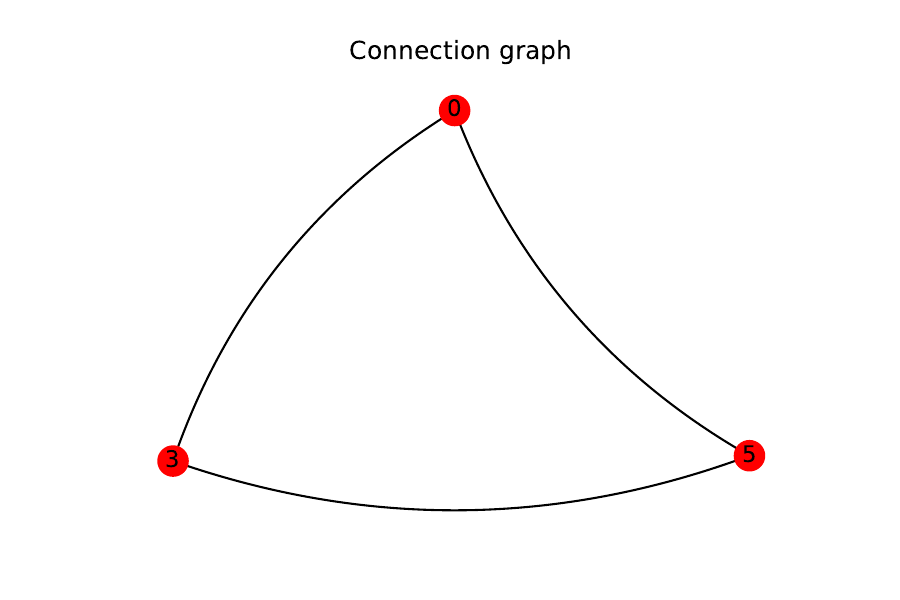}
        \caption{Connection graph \textsf{ConnG}}
    \end{subfigure}
    \hfill
\begin{subfigure}{0.48\textwidth}
        \centering
		\includegraphics[scale = 0.6]{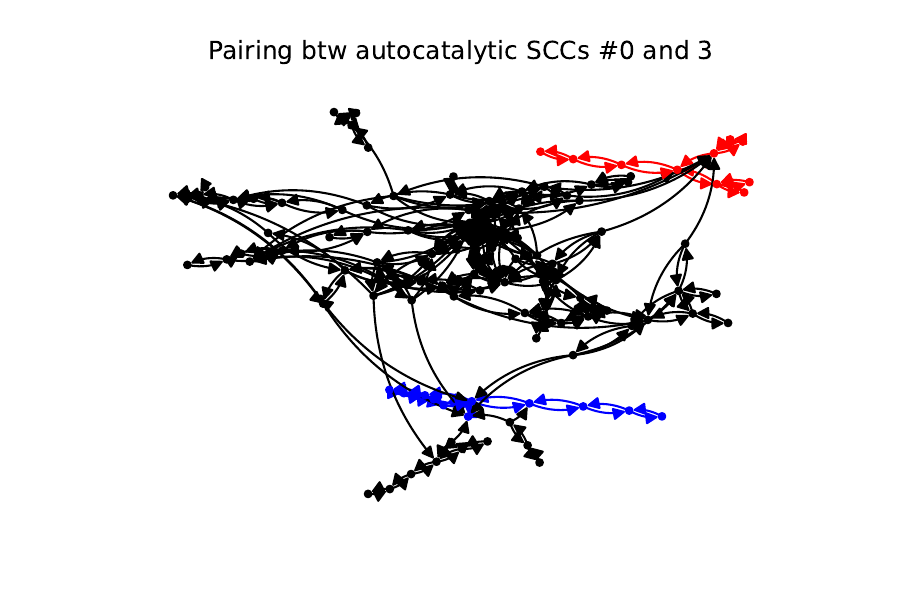}
		\caption{Pairing btw autocatalytic SCCs \#0 and \#3}
	 \end{subfigure}

\caption{Autocatalytic SCCs and their pairings}
\label{random_kinetic_network7_2} 
\end{figure*}

\Medskip We generated ten instances of random networks with $n_{max}=5$. Two of them (resp. 5, 3) contained 
one  (resp. 2, 3) autocatalytic SCC (generating more instances shows that it is quite possible to get 
none). Fig. \ref{random_kinetic_network7_2} reports one instance with 3 autocatalytic 
SCCs indexed 0, 3, 5, see plot (a) for the substrate graph $\cal S$. SCC \#5 is somewhat large (see plot (b)) and would require a larger plot, but SCCs
\#0 and \#3 are small components, exhibiting 3, resp. 4 different addition mechanisms 
(see plots (c) and (d)). All 3 components are connected by \textsf{ConnG} (plot (e)); 
in all ten instances, $\textsf{ConnG}$ is the total graph. 
Plot (f) gives a streamlined plot of their location within ${\cal S}$, with 
${\cal C}_0$ (resp. ${\cal C}_3$)  nodes and edges colored red (resp. blue).

%%%%%%%%%%%%%%%%%%%%%%%%%%%%%%%
%%%%%%%%%%%%%%%%%%%%%%%%%%%%%

\section{Explorations of high-level mechanisms: theory and simulations} \label{section:high-level}

%%%%%%%%%%%%%%%%%%%%%%
%%%%%%%%%%%%%%%%%%%%%%%%

By contrast with \S \ref{section:low-level},
we explore in this section the high-level random reactivity structure of RRBCN models.
This structure is determined by the number and organization of high-level 
mechanisms (equivalently, complex primitive rules). Among the foremost questions, formulated already in the introduction:
are there arbitrarily complex mechanisms, suggesting possible open-ended evolution ? 
if not, what are the most complex mechanisms ? how complex is a 'typical' mechanism ?  
what is the typical shape of a large-level reaction subnetwork ?

\Medskip Simulations are of limited value for such questions. The main tool here has
been an approximate mapping to models of random tree growth, called Galton-Watson models, which
are classical in the mathematics literature (\S \ref{subsection:mapping}). This 
mapping has made it possible to delve into the above questions. We report here mainly
on asymptotic rule statistics (\S \ref{subsection:I} for Model I, \S \ref{subsection:II}
for Model II) which answer the above questions regarding the complexity of mechanisms
within the limits of RRBCN models. Models I and II behave very differently in this
respect: whereas the complexity of mechanisms is capped in Model I, there are regimes
in the parameter space of Model II for which one finds arbitrarily complex mechanisms. 
Regimes present themselves in the form of a phase diagram. This fundamental distinction explains why it has been possible to validate our predictions for Model I by detailed
simulations, while (despite no contradiction) numerical evidence is more elusive for Model II.

\Medskip A large part of our work has been relegated to Suppl. Info., \S II-V, notably
simulations (\S V), with only
main arguments and plots presented in Main Text. Also, work in progress on the typical shape of a large-level reaction subnetwork is presented in Suppl. Info. \S III.B.
(anabolism) and IV.B. (catabolism);
superposing the anabolic and catabolic parts to study the typical shape of reaction subnetworks
is a non-trivial, yet unsolved problem, left for further investigation.

%%%%%%%%%%%%%%%%%%%%%%%%%%%%%%%%%

\subsection{Mapping to tree models}  \label{subsection:mapping}

%%%%%%%%%%%%%%%%%%%%%%%%%%%%%

As pointed out earlier,  parameter values $p\simeq 1$ or $z\simeq 0$ make the large-level random structure
very poor. This suggests we look at intermediate regimes where $p$ is not close to 1 
(typically, $p<0.5$; possibly, $p\to 0$) and $z$ is not too small (typically, $z>0.5$). These conditions
favor the emergence of a small number of unexpectedly complex mechanisms $\rho$
(with $m_{\rho}\ge 10$, say), an interesting
scenario for evolution. 
They also make simulation tools less practical, or sometimes prohibitively time- and 
memory-consuming (the program crashes), because the number of molecules, roughly
 $|{\cal A}|^{n+1}/2$, increases exponentially with  
level $n$. Limiting to exhaustion of primitive rules (see \ref{subsection:RRBCN-I-II}) is necessary, but
not sufficient when the number of those becomes large.

\Medskip 
Fortunately, many large-level observables of RRBCN models may be computed analytically through an  approximately mapping to a random tree model, that we set out to  describe. Anabolism and catabolism can be studied separately and with similar techniques. We deal here with anabolism, fix $f\in {\cal F}$ and 
consider admissible reactions $\rho : f+x\to f\cdot x$. {\em Context completion here is bilateral} for a better general agreement with real-world chemistry: if $f$ reacts with $x$, then it reacts with $r_1\cdot x\cdot r_2$ (and the 
product is branching if $r_1$ is not trivial).  In particular, it reacts with $x\cdot r_2$. Representing the set of all molecules $x=A_1 A_2\cdots$ as a tree, see Fig.
\ref{fig:max-compo-tree} for ${\cal A} = \{A,B,C\}$,  keeping only those $x$ that
do {\em not} react with $f$ (in red), and those such that $\rho$ is a {\em primitive} rule (in blue), yields a {\em composition tree}. By construction, a composition tree
is a tree $\T=\T(f)$ with (formal) root $\emptyset$, red inner vertices and blue leaves. 
The maximal composition tree (containing all words, in red) is obtained when $f$ does not react
with any $x$.          
 
\bigskip
\begin{figure}[H] 
\begin{center}
\begin{tikzpicture}[scale = 0.6]   

\draw(6,4) node {$|x|=0$};
\draw(6,4+4*1/2) node {$|x|=1$};

%\draw(0,0) node {\textbullet}; 
\draw(0,4) node {\textbullet}; \draw(0.5,4) node {$B$};
\draw(0,4+4*1/2) node {\textbullet}; \draw(0,4+4*1/2+0.5) node {\small $BB$};

\draw(0,0)--(0,4)--(0,4+4*1/2); 
\draw[dashed](0,4+4*1/2)--(0,4+4*1/2+4*1/4);
\draw(0,0)--(-4*0.707,4); \draw(-4*0.707,4) node {\textbullet};
\draw(-0.5-4*0.707,4) node {$A$};
\draw(0,0)--(4*0.707,4);
\draw(4*0.707,4) node {\textbullet};
\draw(0.5+4*0.707,4) node {$C$};

%%%%%%%%%%%%%%%%%%%%%%
\draw(-4*0.707,4)--(-4*0.707,4+4*1/2);
\draw(-4*0.707,4)--(-4*0.707-4*0.707*1/3,4+4*1/2);
\draw(-4*0.707,4)--(-4*0.707+4*0.707*1/3,4+4*1/2);

\draw(-4*0.707,4+4*1/2) node {\textbullet};
\draw(-4*0.707+4*0.707*1/3,4+4*1/2) node {\textbullet};
\draw(-4*0.707-4*0.707*1/3,4+4*1/2) node {\textbullet};

\draw(-4*0.707,4+4*1/2+0.5) node{\small $AB$};
\draw(-0.5-4*0.707-4*0.707*1/3,4+4*1/2) node {\small $AA$};
\draw(-4*0.707+4*0.707*1/3,4+4*1/2+0.5) node {\small $AC$};
\draw(0,4)--(-4*0.707*1/3,4+4*0.5);
\draw(-4*0.707*1/3,4+4*0.5) node {\textbullet};
\draw(0,4)--(4*0.707*1/3,4+4*0.5);
\draw(4*0.707*1/3,4+4*0.5) node {\textbullet};
\draw(-4*0.707*1/3,4+4*0.5+0.5) node {\small $BA$};
\draw(4*0.707*1/3,4+4*0.5+0.5) node {\small $BC$};

\draw(4*0.707,4)--(4*0.707,4+4*1/2);
\draw(4*0.707,4)--(4*0.707-4*0.707*1/3,4+4*1/2);
\draw(4*0.707,4)--(4*0.707+4*0.707*1/3,4+4*1/2);
\draw(4*0.707,4+4*1/2) node {\textbullet};
\draw(4*0.707+4*0.707*1/3,4+4*1/2) node {\textbullet};
\draw(4*0.707-4*0.707*1/3,4+4*1/2) node {\textbullet};
\draw(4*0.707,4+4*1/2+0.5) node{\small $CB$};
\draw(0.5+4*0.707+4*0.707*1/3,4+4*1/2) node {\small $CC$};
\draw(4*0.707-4*0.707*1/3,4+4*1/2+0.5) node {\small $CA$};

%%%%%%%%%%%%%%%%%%%%%%%%%
\draw[dashed](-4*0.707-4*0.707*1/3,4+4*1/2)--(-4*0.707-4*0.707*1/3-4*0.707*1/9, 4+4*1/2+4*1/4);

\draw[dashed](-4*0.707,4+4*1/2)--(-4*0.707,4+4*1/2+4*1/4);

\draw[dashed](-4*0.707+4*0.707*1/3,4+4*1/2)--(-4*0.707+4*0.707*1/3+4*0.707*1/9, 4+4*1/2+4*1/4);

\draw[dashed](4*0.707+4*0.707*1/3,4+4*1/2)--(4*0.707+4*0.707*1/3+4*0.707*1/9, 4+4*1/2+4*1/4);

\draw[dashed](4*0.707,4+4*1/2)--(4*0.707,4+4*1/2+4*1/4);

\draw[dashed](4*0.707-4*0.707*1/3,4+4*1/2)--(4*0.707-4*0.707*1/3-4*0.707*1/9, 4+4*1/2+4*1/4);

\end{tikzpicture}
\end{center}

\caption{Maximal  composition tree $\T$, first levels: $V_0(\T) = \{A,B,C\},\ 
V_1(\T) = \{AA,AB,AC,BA,BB,BC,CA,CB,CC\}.$}
\label{fig:max-compo-tree}
\end{figure}

  Equivalently, letting $\omega_x^f$ be the decision variable attached to the rule
  $f+x\to f\cdot x$:
\BEQ (x\in \T(f)) \Leftrightarrow \Big(\forall 
x'\subsetneq x, \omega^f_{x'} = 0\Big) \label{eq:xTf}
\EEQ
and $x\in \T(f)$ is a leaf iff $\omega^f_x=1$.  See Suppl. Info. \S IV. A. for  a similar
construction applied to catabolism in terms of {\em fragmentation tress}.

\Medskip Our predictions, consistent with available (partial but still very substantial) simulations,
are built on an  approximate mapping from our model to an {\em inhomogeneous Galton-Watson (GW) tree model}.  We refer the reader to the book \cite{AthNey}, Chapter 1 \& 5 for a standard introduction
to the subject. GW trees are also used to represent the different generations of cellular division models in biology \cite{KobSug}. By definition, individuals are indexed by an integer generation number. Each generation $n$ individual $(n,j)$ has a random number of descendants $\zeta_{(n,j)}$ according to a fixed {\em progeny law}; the model is called
inhomogeneous if that law depends on $n$. All progenies $\zeta_{(n,j)}$ are assumed to be independent. Moments of the model may be computed using generating
functions $f_n(s) = \esper[s^{\zeta_n}]$, $s\ge 0$, where $\zeta_n$ is any of the
$\zeta_{(n,j)}$ variables.

\Medskip Our mapping is based  on a 
$(1/|{\cal A}|)$ expansion for $|{\cal A}|$ large. We denote by ${\cal X}_n$ the 
set of level $n$ words, and call (primitive) $f$-reactants those $x$ such that 
$\rho : f+x\to f\cdot x$ is a (primitive) rule. We observe the following:

\Bigskip {\em Key observation 1.} Let $x=A_1\cdots A_n\in {\cal X}_{n-1}$, and $A\in {\cal A}$. 
Assume $x\in \T(f)$. {\em Then $x\cdot A\in \T(f)$ if and only if for all $k=0,\ldots,n$,  
$A_{n-k+1}\ldots A_n\cdot A$ is not a primitive $f$-reactant.} (By abuse of notation,  we let $A_{n-k+1}\ldots A_n \cdot A = A$  
if $k=0$).  Namely, as already observed, all other subwords of $x\cdot A$ are actually subwords of $w$,
therefore are not $f$-reactants since $x\in \T(f)$.

\Medskip {\em Key observation 2.} Choose $w=A_1\cdots A_n$ uniformly in 
\BEQ V_{n-1}(\T) = \{x\in \T \ |\ |x|=n-1\}.
\EEQ
 Denote
$\zeta_x = (\zeta_x^A)_{a\in {\cal G}}$ its progeny; in other words, $\zeta_x^A=1$ iff 
$x\cdot A\in \T$. 
 Then, {\em to leading order in $1/|{\cal A}|$, the progeny law of $x$ has independent components, and the law of each component is independent of the choice of $x$.}  Namely, the $(n+1)\times|{\cal A}|$ variables 
 \BEQ \omega^f_k(x,A) \equiv \omega^f_{A_{n-k+1}\ldots A_n\cdot A}, \qquad 0\le k\le n, \ 
A\in {\cal A}
\EEQ
 (see \S \ref{eq:xTf}) are independent since they are indexed by different words. Furthermore, $\omega^f_{A_{n-k+1}\ldots A_n\cdot A}$ depends on $w$ only through the possibility that $A_{n-k+1}\ldots A_n\cdot A$ is a subword of $w$. However, fixing these $k+1$ atoms implies restraining to a proportion $\approx |{\cal A}|^{-(k+1)}$ of
$V_n(\T)$. For  $n=O( |{\cal A}|/p)$, this implies (considering all potential $n-k$ 
locations of $A_{n-k+1}\ldots A_n\cdot A$ in $w$, and taking into account the variance $p(1-p)$ of a Bernoulli law) a correcting term to the progeny of order $p(1-p)\times   |{\cal A}|^{-(k+1)} 
\times n = O(|{\cal A}|^{-1})$,  {\em except} for $k=0$.   The correction associated to the 
$(k=0)$ term can be computed, and -- as turns out -- does not alter
significantly the consequences of the above observations.

\Medskip The mapping to a GW model resulting from these key observations is described
in \S \ref{subsection:I} for Model I, and \S  \ref{subsection:II} for Model II.
The approximate mapping makes statistical quantities for the large-level random reactivity structure of RRBCN 
models analytically computable. We present the results in the Main Text, and refer the
reader to Suppl. Info. for proofs. 

\Medskip The mapping turns exact when $|{\cal A}|\to\infty$. This makes it particularly
well-suited  for  applications to chemistry involving some coarse-graining, since $\cal A$ should span 
a set of functional groups of size 20-30 at least. Actually, simulations show a surprisingly good agreement with our estimates in the case of Model I (see Suppl. Info. 
\S V.A. and \S V.C.) even when $|{\cal A}|=3$, describing e.g. heavy atoms for
CHON molecules.

%%%%%%%%%%%%%%%%%%%%%%%ù
%%%%%%%%%%%%%%%%%%%%%%

%%%%%%%%%%%%%%%%%%%%%%%%%

\subsection{Asymptotic rule statistics: the case of Model I}
\label{subsection:I}

%%%%%%%%%%%%%%%%%%%%%%%%%

We discuss here the structure of  large-level anabolic primitive rules for Model I, as deduced from  asymptotics for anabolic composition trees (\S \ref{subsubsection:I-ana}).
Details and proofs may be found in Suppl. Info. \S II. Results for catabolis are
briefly discussed in \S \ref{subsubsection:I-cata}, see Suppl. Info \S IV.

%%%%%%%%%%%%%%%%%%%%%%       
       
\subsubsection{Anabolism} \label{subsubsection:I-ana}      
       
%%%%%%%%%%%%%%%%%%%%%%       
       
In the whole paragraph, we fix a composition tree 
$\T=\T(f)$.  
Dealing simultaneously with several foodset molecules is a very easy generalization; see 
next section (Model II), where we start directly from a general foodset $\cal F$.

\Medskip  Taking into account the two Key observations, and forgetting $(k=0)$ corrections,  one can expect that, for generic $w\in V_{n-1}(\T)$,  $w.A\in  V_{n}(\T)$ with probability $\sim \proba[\omega^f_k(x,A)=0, \ 1\le k\le n] = (1-p)^n$. Thus, the 
mapping of Model I is to the inhomogeneous GW model with generation $n$ progeny law equal
to the binomial law
\BEQ (progeny\ law\ of\ GW\ model) \qquad Bin(|{\cal A}|,(1-p)^n)
\EEQ

\Bigskip In more concrete terms,  
 {\em the average number of
descendants of a typical vertex of level $n-1$ in a composition
tree is roughly $ \sim |{\cal A}| (1-p)^n.$}  
This is to be understood in a weak sense (see below), namely, 
$\log \langle |V_n(\T)|\rangle \sim \log \prod_{k=1}^{n} (|{\cal A}|(1-p)^k)$
for $n$ large enough but still less than $O(|{\cal A}|/p)$, or
\BEQ \log \langle |V_n(\T)|\rangle \sim  n \Big( \log  |{\cal A}| - \frac{n+1}{2} \log(1/(1-p)) \Big)
  .  \label{eq:I-logVn}
\EEQ 
       
\Medskip Thus composition trees dwindle rapidly, and become rapidly
extinct as $n$ exceeds the value of $n$ which maximizes (\ref{eq:I-logVn}), 
\BEQ   {\mathrm{(maximum\ average\ size\ level) }}
 \qquad n_0:= \frac{\log(|{\cal A}|)}{\log(1/(1-p))}, \label{eq:n0}
\EEQ
solution of the equation 
\BEQ (1-p)^{n_0} |{\cal A}| = 1.
\EEQ
In practice, since we are interested in small values of  $p$, we often use the approximate value
 \BEQ n_0\sim \frac{\log(|{\cal A}|)}{p}. \EEQ 
Note that the r.-h.s. in (\ref{eq:I-logVn}) becomes negative when $n>n'_0$, with
\BEQ {\mathrm{(extinction\ level)}} \qquad  n'_0:= 2n_0 \label{eq:n'0}
 \EEQ

\Bigskip The actual interest is not in the composition tree, but in its leaves, which 
are the primitive $f$-reactants. We characterize the set of leaves
\BEQ {\mathrm{Prim}} = \cup_{n\ge 0} {\mathrm{Prim}}_n, \qquad 
 {\mathrm{Prim}}_n = {\mathrm{Prim}}  \cap V_n(\T) \label{eq:Prim_n}
\EEQ
 by its shape: size (number of primitive rules) by level, and top level, which is the maximum complexity index of all $f$-addition 
mechanisms.

\begin{figure}[t]
\centering
\includegraphics[scale = 0.6]{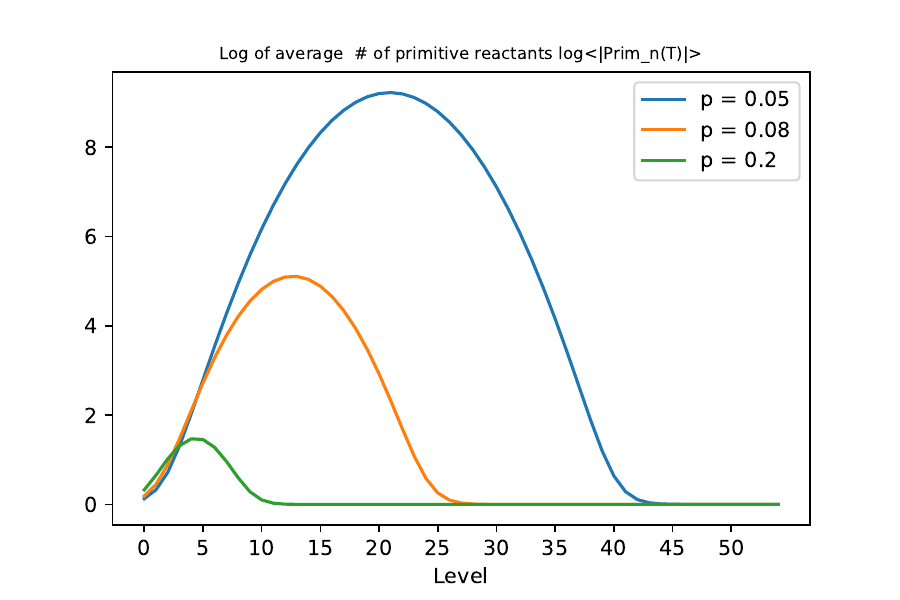}  
\caption{Log of average number of primitive reactants in anabolic Model I for 
$|{\cal A}|=3$ }
\label{fig:log-prim-I}
\end{figure}

\Bigskip {\bf Number of primitive rules by level} (see Suppl. Info., \S II.A. and 
\S II.B.).  Let $x\in V_{n-1}(\T)$. Recall (see Key
 observation 2) that $f+x\cdot A\to f\cdot x\cdot A$ is a primitive reaction iff $\omega^f_k(x,A) = 0$, $1\le k\le n-1$, and $\omega^f_n(x,A) \equiv \omega^f_{A_1\cdots A_n\cdot A } = 1$.
 Thus the number of blue vertices of level $n$ may be computed as above, except that
the last generation progeny law is not Bin$(|{\cal A}|,p_n)$ but Bin$(|{\cal A}|,p(1-p)^{n-1})$. Taking
into account the constant prefactor $|{\cal A}|(1-p)^{|{\cal A}|}$ coming from the $(k=0)$ correction (see Suppl. Info \S II.B.), this 
implies:
\BEA \langle |\Prim_n|\rangle &\approx & |{\cal A}| (1-p)^{|{\cal A}|} \, \times \prod_{k=1}^{n-1} (|{\cal A}|
(1-p)^k)  \nonumber\\
 && \qquad \qquad \times (|{\cal A}| p(1-p)^{n-1})  \nonumber\\
&=&  p |{\cal A}|^{n+1} (1-p)^{|{\cal A}|-1 + n(n+1)/2} 
\label{eq:Primn}
\EEA
Just like $\langle V_n(\T)\rangle$, it is maximum for $n\simeq n_0$, see (\ref{eq:n0}), 
which can therefore be thought as the 'typical' complexity of a specific anabolic mechanism. Figure \ref{fig:log-prim-I} gives the number of primitive rules in terms of $n$ for various values of $p$.

\Medskip The average total number of primitive reactions is $\langle |\Prim| \rangle= \sum_{n\ge 0} \langle|\Prim_n|\rangle$. The largest term in the summand is obtained for $n\simeq n_0$, yielding 
$\langle |\Prim|\rangle \approx p|{\cal A}|^{n_0+1}  (1-p)^{|{\cal A}|-1 + n_0(n_0+1)/2} 
\sim p |{\cal A}|^{1+ \log|{\cal A}|/2p}.$

\Bigskip  {\bf Top level and extinction probabilities} (see Suppl. Info., \S II.C. and \S II.D.). The {\em top level} (or {\em height} of $\T$, denoted ht$(\T)$) is the level at which the tree $\T$ becomes
extinct. Let $(Z_n)_{n\ge 0}$ be the number of generation $n$ individuals in the GW model.    Denote by $P(m,j;n,i):=\proba[Z_n =i \ |\ Z_m =j ]$  $(n>m)$ the inhomogeneous iterated  
kernel associated to the Markov process $(Z_n)_{n\ge 0}$. By definition, $\T$ is extinct
 at time $n$ with probability $u_{0,n}:=P(0,1;n,0)$. The Chapman-Kolmogorov equation implies
 $u_{0,n} = \sum_{k\ge 0}  P(0,1;1,k) P(1,k;n-1,0)$; then $P(1,k;n-1,0) = (P(1,1;n-1,0))^k$ by 
 independence of the branches. Denote by $u_{k,n} := P(k,1;n,0)$ the probability to be extinct at time $n$
 starting from one individual at time $k$; then we have found $u_{0,n} = f_1(u_{1,n})$. Note that $u_{n-1,n} = (1-p_n)^{|{\cal A}|} = f_n(0)$.  Iterating, we
 get $u_{0,n} = f^{(n-1)}(u_{n-1,n}) = f^{(n)}(0)$; thus, 
\BEQ u_{0,n} = f^{(n-1)}((1-p_n)^{|{\cal A}|}).  \label{eq:u0n} 
\EEQ

Extinction probabilities are estimated as a by-product of a general estimation for the iterated 
function $f^{(n)}(s), 0<s<1$, see Suppl. Info. \S II.  B., where it is proves that trees grow extinct with high probability at a level $\sim 2n_0$. Thus the top level (maximum mechanism complexity) is roughly equal to the extinction level $n'_0$ of 
eq. (\ref{eq:n'0}).

\Bigskip {\bf Logarithmic size} (see Suppl. Info., \S II.E.).  One may be interested in measuring 
and estimating $\langle \log |V_n(\T)| \rangle$ instead of $\log \langle |V_n(\T)|\rangle$, and $\langle \log (|\Prim_n|)\rangle$ instead of $\log \langle |\Prim_n|\rangle$. This induces in general a supplementary prefactor $<1$, since (by convexity) $\langle \log(X)\rangle \le \log \langle X\rangle$ for any
positive random variable $X$. This is a well-known issue in the study of disordered media; interpreting
$|V_n(\T)|$  as a random partition function, $\langle  \log |V_n(\T)| \rangle \sim \log \langle |V_n(\T)|\rangle$ for large $n$ is indicative of a self-averaging measure, i.e. the distribution of $|V_n(\T)|$ is very concentrated. True asymptotics $n\to\infty$ are however irrelevant here, since trees go extinct.

\Medskip It is proved in  Suppl. Info., \S II. C., that self-averaging holds for 
the Galton-Watson model, but not for our original  model.

%%%%%%%%%%%%%%%%%%%%%%       
       
\subsubsection{Catabolism} \label{subsubsection:I-cata}      
       
%%%%%%%%%%%%%%%%%%%%%%   

We only sketch our results, and refer the reader to Suppl. Info. \S IV.A. for
more information. Instead of a set of $f$-dependent composition trees $\T(f)$, we get
a unique {\em fragmentation tree} $\T_{frag}$, whose leaves coincide with the set 
of $x$ such that $x$ is decomposable by a primitive rule, i.e. at least one catabolic reaction of the
form $x\to x_1+ x_2$ ($x_1\cdot x_2= x$) is a primitive rule of the RRBCN-network $\bar{\cal R}$. Then maximum average size level $n_0$  (\ref{eq:n0}), resp. extinction level $n'_0$
(\ref{eq:n'0}) have equivalent counterparts $n_{frag}$, resp. $n'_{frag}$ for catabolism, see Suppl. Info., eq. (IV.6), 
\BEQ n_{frag} \sim \sqrt{\frac{2\log|{\cal A|}}{q}}, \qquad n'_{frag}\sim \sqrt{3}\ n_{frag}
\EEQ
Multiple fragmentation processes a priori available for any given $x$ explain why 
the complexity of catabolic mechanisms is typically much lower than in the anabolic case, as evidenced by the square-root in $n_{frag}$.

%%%%%%%%%%%%%%%%%%%%%%%%%%%%%

\subsection{Asymptotic rule statistics : the case of Model II}
\label{subsection:II}

%%%%%%%%%%%%%%%%%%%%%%%%%

\begin{figure}[t]
\centering
\includegraphics[scale = 0.6]{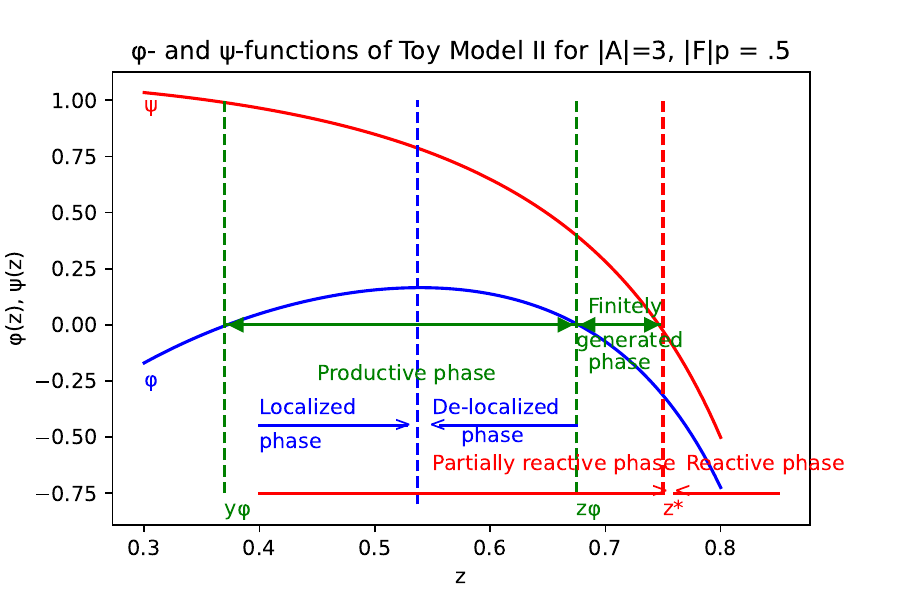}  
\caption{Anabolic phases of Model II for 
$|{\cal A}|=3$ and $|{\cal F}|p = 0.5$}
\label{fig:phases}
\end{figure}

We now deal with Model II by reusing the same techniques. The analysis is more involved than for Model I. Our main 
finding is a phase diagram in the parameter set $(p, q, z, |{\cal A}|, |{\cal F}|)$  for the structure of large-scale primitive rules, see Fig. \ref{fig:phases}, including
a phase in which the number of primitive rules is infinite, which hints at the eventuality for open-ended evolution. We discuss the case of anabolism   in some detail, but refer to Suppl. Info. \S III for
complements and proofs. The case of catabolism, which is similar, is left out entirely,  see Suppl. Info. \S IV.A, in particular, the catabolic phase diagram, Fig. IV. 13.

\Medskip The fugacity parameter $z$ is crucial here. In particular, 
for $z$ not too large, the composition tree is infinite; more precisely, the average fraction
 of molecules of level $n$ that react is asymptotically of the form $(c(z))^n$ for 
 some constant $c(z)\in (0,1)$. This makes Model II look more  realistic. Yet another distinctive feature of this model is that it admits an
infinite number of primitive rules in an intermediate regime of values of $z$. This phase, particulary appealing for chemical evolution, is
called {\bf productive phase}. Statistical features within each phase are conveyed by
two explicit {\bf phase functions} $\psi_{ana},\phi_{ana}$.

\Medskip We consider here all elements in the foodset and derive properties of a global composition
 tree $\T$. 
 
\Medskip One significant difference between Model I and Model II is the existence of thresholds depending
on the parameter $z$. This may be seen immediately by looking at the probability that
$\bar{\cal R}$ is empty (no reaction at all), formally,
\BEQ \proba[\bar{\cal R} = \emptyset] = \proba[\forall x\in {\cal X},\forall f\in {\cal F},\ \omega^f_x=0] = 
\prod_{n\ge 0} (1-pz^n)^{|{\cal F}|\, |{\cal A}|^{n+1}} \EEQ
Then 
\BEQ -\log(\proba[\bar{\cal R}(\omega) = \emptyset]) \sim p|{\cal A}||\, |{\cal F}|\, \sum_{n\ge 0} 
|{\cal A}|^n z^n 
\EEQ
The series diverges if and only if $z\ge 1/|{\cal A}|$. Thus the above probability is 0 above the
threshold value $z=1/|{\cal A}|$, including the case of Model I $(z=1)$. By contrast, if 
$z<1/|{\cal A}|$, there is a nonzero probability that the reaction network is empty. For those values of $z$, it may be interesting
to consider the measure conditioned to the event $\bar{\cal R}(\omega) \not= \emptyset$ and see how the
results below change. In practice, we tend to favor in the following analysis (at least in simulations)  values of $z$ close to 1.

\Bigskip Let us reproduce the key arguments of \S \ref{subsection:I}, taking
into account all foodset elements this time. We let $\T:= \cap_{F\in {\cal F}}\T(F)$; this  composition tree may be constructed
 in the following way:   if $x=A_1\cdots A_n\in V_{n-1}(\T)$, then $x\cdot A\in V_n(\T)$
if and only if, for every $F\in {\cal F}$, the reaction $f + x\cdot A \to f\cdot x\cdot A$ is not in $\bar{\cal R}$. 
Key observations 1. and 2. still hold, with simple modifications; for a generic $x\in V_{n-1}(\T)$, $x\cdot A 
\in V_n(\T)$ with probability 
\BEA p_n &\sim&  \prod_{k=1}^n \Big\{\prod_{F\in {\cal F}} (1-p_{f +A_{n-k+1}\cdots A_n\cdot A
\to f\cdot A_{n-k+1}\cdots A_n\cdot A}) \Big\}  \nonumber\\
&=& \prod_{k=1}^n \prod_{f\in {\cal F}}  (1-p'_f z^k)
\label{eq:II-pn},
\EEA
 where $p'_f= p'_f(z):= pz^{|f|+1}$.  This allows a comparison to an inhomogeneous Galton-Watson tree  with generation $n$ 
law Ber($|{\cal A}|, p_n)$, where $p_n:= \prod_{k=1}^n  \prod_{f\in {\cal F}} (1-p'_f z^k)$ instead of 
$(1-p)^n$. Taking into account the $(k=0)$ correction, we get a multiplicative factor $|{\cal A}|$ (maximal number of branches), times $ \prod_{f\in {\cal F}} (1-p'_f)^{|{\cal A}|}$ (probability for all reactions $f+A\to 
f\cdot A$, \ $(f,A)\in {\cal F}\times {\cal A}$ to be rejected).  We denote by $Z_n$ the size of the $n$-th generation.
 Hence our prediction is
\BEA && \langle |V_n(\T)|\rangle \approx  \esper[Z_n] =  |{\cal A}| \ \prod_{f\in {\cal F}} 
(1-p'_f)^{|{\cal A}|} \ \cdot\ \prod_{j=1}^n (|{\cal A}|p_j)  \nonumber\\
&& \ \ =   |{\cal A}|^{n+1} \ 
\prod_{f\in {\cal F}} (1-p'_f)^{|{\cal A}|} \ \times\ \prod_{k=1}^{n}\prod_{f\in {\cal F}}
 (1- p'_f z^k)^{(n-k+1)} \nonumber\\  \label{eq:II-Vn}
\EEA
hence (when $p'\to 0$ and $n$ is large), letting $p'=p'(z):= p \langle z^{|F|+1}\rangle_{\cal F} = 
\frac{p}{|{\cal F}|} \sum_{f\in {\cal F}} z^{|f|+1}$, 

\BEA && \log \langle |V_n(\T)|\rangle \sim -p'|{\cal F}||{\cal A}|   \nonumber\\
&& \ \  +  n \Big\{ \log  |{\cal A}| - |{\cal F}| p'\sum_{k=1}^{n} (1-\frac{k-1}{n}) 
z^k  \Big\}  \label{eq:II-logVn}
\EEA 

\Medskip 
The main difference with the case of Model I is that, when $z<1$,   the sum in $k$ replacing
the term $\frac{n+1}{2}$ in the r.-h.s. of  (\ref{eq:I-logVn}) is bounded in $n$.
Assuming  $z<c<1$ to be bounded away from 1, and  
replacing each term in the above sum in $k$ by its equivalent when $n\to\infty$, yields
 \BEA && \log \langle |V_n(\T)|\rangle \sim  -p'|{\cal F}||{\cal A}|  +  n \psi_{ana}(z), \nonumber\\
&&\qquad   \psi_{ana}(z)\equiv \psi_{ana}(z|\, |{\cal A}|, |{\cal F}|, p') := \log|{\cal A}| - |{\cal F}| p' \frac{z}{1-z} \nonumber\\ 
\EEA
Thus, when $\psi_{ana}(z)>0$,  the average fraction
 of molecules of level $n$ that react is $\approx c(z)^n$,  with $\log(c(z)) = - |{\cal F}| p' \frac{z}{1-z}$. 

%%%%%%%%%%%%%%%%%%%%%%%%%

\Medskip {\em Primitive vertices.} Similarly, $\langle |\Prim_n(\T)| \rangle$ may be estimated by 
replacing the factor $\prod_{f\in {\cal F}} (1-p'_f z^n)$ in the expression (\ref{eq:II-pn}) for $p_n$ with
$1- \prod_{f\in {\cal F}} (1-p'_f z^n) \sim |{\cal F}|p' z^n$.   Thus

\BEA &&  \langle |\Prim_n|\rangle \approx |{\cal F}| p' z^n\  |{\cal A}|^{n+1}
\prod_{f\in {\cal F}}  (1-p'_f)^{|{\cal A}|}
\, \times\nonumber\\
&&\qquad \prod_{k=1}^{n-1} \prod_{f\in {\cal F}}
 (1-p'_f z^k)^{(n-k+1)}   \label{eq:II-Primn}
\EEA

whence

\BEA && \log \langle |\Prim_n|\rangle \sim  -p'|{\cal F}||{\cal A}| + \log(|{\cal F}| p') + n \phi_{ana}(z),
\nonumber\\ && \qquad \phi_{ana}(z):= \log(z)+ 
\psi_{ana}(z)
\EEA

\Medskip The {\bf anabolic phase functions} $\psi_{ana}=\psi_{ana}(z), \phi_{ana}=\phi_{ana}(z)$ are studied in Suppl. Info., \S III in the particular case when $p'(z)=pz$ (i.e. if 
$|F|=0$ for all $F\in {\cal F}$), a convenient choice for computations; other choices lead 
to very similar conclusions.

\Bigskip Letting $z$ increase from $0$ to $1$, we potentially  get three different phases (depending on the  parameters, they may
or may not be all existent):

\Medskip {\bf (a) (Partially reactive, productive phase).} This phase is  defined by the set of parameters $(p,z)$ such that   $\phi_{ana}(z)>0$ (implying, in particular, $\psi_{ana}(z)>0$), or:
\BEQ \frac{p' z}{1-z} < \frac{\log(|{\cal A}| z)}{|{\cal F}|} \EEQ 
It is 
called  {\em productive},  since one finds new primitive reactions at every level, though the number of these
grows much slower than the number of words. (A) is called: {\em partially} reactive, because not all
molecules of high enough level react. Assume $p'(z)=pz$.  If  $\max_{z\in (0,1])}\phi_{ana}(z)\le 0$, in particular, if 
$|{\cal F}|p>\frac{|{\cal A}|^2}{2e}$ (see Lemma II. 1 (iv) in Suppl. Info. \S III),
 then this phase is empty;
however, if $|{\cal F}|p\ll 1$, then, by (ii) of the same Lemma  (see also (iii) for a more quantitative
criterion), phase (A) is realized for $0<y_{\phi}<z<z_{\phi}<1$ where $y_{\phi},z_{\phi}$ are the two zeros of $\phi$ on $(0,1)$, 
with $1-z_{\phi} \sim_{p\to 0} \frac{|{\cal F}|p}{\log|{\cal A}|}$. Under the same condition, phase (A) splits into two qualitatively different subphases:  

\textbullet\ a {\em localized phase} defined by $z<z'_{\phi}\sim_{p\to 0} 1-\sqrt{|{\cal F}|p}$ and characterized by $\phi'_{ana}(z)>0$, which thus behaves similarly to the free random model (no context derivation, all admissible
reactions are primitive) for which $\phi_{ana}(z) = \log(z)$;

\textbullet\ a {\em (delocalized phase)}, when $z'_{\phi}<z<z_{\phi}$, $\phi'(z)<0$; this indicates a very different subphase,
in which there are  less and less high level primitive rules as $z$ increases, because more and more rules are derived from low-level primitive rules. At $z=z'_{\phi}$, we get $\phi_{ana}(z'_{\phi})\sim \log|{\cal A}| - 2
\sqrt{|{\cal F}|p}$, so that the fraction of primitive reactants of level $n$ is  $\frac{\langle |\Prim_n|\rangle}{|{\cal A}|^n} \approx (e^{-2\sqrt{|{\cal F}|p}})^n$.

%%%%%%%%%%%%%

\Medskip {\bf  (b) (Partially reactive, finitely generated phase).} This phase is  defined by the set of parameters $(p,z)$ such that   $\phi_{ana}(z)<0$ but   $\psi_{ana}(z)>0$, or:
\BEQ  \frac{\log(|{\cal A}| z)}{|{\cal F}|} < \frac{p' z}{1-z} < \frac{\log|{\cal A}| }{|{\cal F}|} 
\EEQ
If (A) is not empty, phase (B) is defined by $z\in (0,y_{\phi})\uplus (z_{\phi},z^*)$ for some $z^* \in (z_{\phi},1)$; it is the disjoint union of a localized part $(z<y_{\phi})$ with
a delocalized one $(z_{\phi}<z<z^*)$.
 Phase (B) is characterized by an exponential decrease
of $\langle |\Prim_n|\rangle$ when $n$ is large, implying that the number of primitive rules
is a.s. finite. As in (A), the set of primitive rules is not large enough to allow all molecules of 
high enough level to react.

%%%%%%%%%%%%%%%%%%
 
\Medskip {\bf (c) (Reactive phase).}   This phase is  defined by the set of parameters $(p,z)$ such that $\psi_{ana}(z)<0$ (implying: $\phi_{ana}(z)<0$), or:
\BEQ   \frac{p' z}{1-z} > \frac{\log|{\cal A}| }{|{\cal F}|}   \label{eq:cond511}
\EEQ
It is similar to Model I, in the sense that there is a finite number of primitive rules, which is large
enough to allow all molecules of 
high enough level to react. We may define a {\em maximum average size level} $n_0$ as in Model I by
requiring that  $ \log \frac{|V_{n+1}(\T)|}{|V_n(\T)|} \sim  \frac{\del }{\del n} \Big( 
n(\log|{\cal A}| - |{\cal F}|p' z \sum_{k=0}^{n-1} z^k) + |{\cal F}|p' z \sum_{k=0}^{n-1} kz^k
\Big)$ vanishes, where $\frac{\del}{\del n}f(n)\equiv f(n+1)-f(n)$ is the discrete derivative. 
This is equivalent to
$\log|{\cal A}| \sim  |{\cal F}|p' z\,  \sum_{k=0}^{n-1}  z^k = |{\cal F}|p' z \frac{1-z^n}{1-z}$, yielding for $\frac{p'z}{1-z} \gg \frac{\log|{\cal A}| }{|{\cal F}|}$
\BEQ n_0 \sim  \frac{\log|{\cal A}|}{|{\cal F}|p'} \EEQ
Similarly, we define the {\em extinction level} $n'_0$ as the level for which the quantity
in curly brackets in  (\ref{eq:II-logVn})
(which is strictly decreasing) changes sign. Computing, we get 
\BEA \log|{\cal A}| &\sim&  |{\cal F}|p' z \sum_{k=0}^{n-1} (1-\frac{k}{n})z^k = 
|{\cal F}|p'z (1-\frac{1}{n}z\partial_z) (\frac{1-z^n}{1-z})  \nonumber\\
&=& \frac{|{\cal F}|p'z}{1-z} 
\, \times\, \Big\{ 1 - \frac{z}{n} \frac{1-z^n}{1-z}\Big\}.
\EEA
Let $\eps:=1-z\to 0$; under the condition (\ref{eq:cond511}), the quantity between curly 
brackets must be small, which implies $1-z^n = 1-(1-\eps)^n \sim n\eps - \frac{n^2\eps^2}{2} \ll 1$, and  then, $1 - \frac{z}{n} \frac{1-z^n}{1-z} \sim \frac{n\eps}{2}$, whence
\BEQ n'_0 \sim 2n_0  \label{eq:5.14}
\EEQ
exactly as in Model I.

\Bigskip We plot $n_0= n_0(z)$ (in purple) in the reactive phase in  Figure 
\ref{fig:n0} below, together with 
$n'_0=n'_0(z)$ (in red). Both quantities diverge as $z$ converges to the lower boundary 
$z^*$ of the reactive phase
($z^*=0.746...$ here). Note that $n'_0(z)/n_0(z)$ explodes as 
$z\searrow z_{\phi}$, so that (\ref{eq:5.14}) fails close to $z_{\phi}$.  The region in the $z$-parameter space
where $n_0$ is large, say $n_0>10$, is small.

\Bigskip It is also instructive to define levels $m_0=m_0(z), m'_0 =m'_0(z)$ similar to $n_0,n'_0$ for
primitive reactions in phases (B)+(C), i.e. in the interval $z\in (z_{\phi},1)$ where $\phi<0$ 
($z_{\phi}=0.676...$ here), namely: $0\equiv \frac{\log|\Prim_{m_0+1}(\T)|}{\log|\Prim_{m_0}(\T)|}$, yielding
 \BEQ m_0 \sim \frac{\log(|{\cal A}|z)}{|{\cal F}|p'z} \EEQ
 (curve in blue),
and $m'_0=m'_0(z)$ (in green) defined by $\log  (|{\cal A}|z) - |{\cal F}| p'z
\sum_{k=0}^{m'_0-1} (1-\frac{k}{m'_0}) 
z^k=0$.

\Bigskip The relatively high value of $|{\cal F}|p$ chosen in the  above Figure is much too large for the above $p\to 0$ estimates to
be even vaguely correct (in particular,  $\phi(z'_{\phi})\simeq 0.165>0$ but $ \log|{\cal A}| - 2
\sqrt{|{\cal F}|p}<0$). However, this choice has made it possible to plot a large part of the
curves $\phi,\psi$ and observe the transition as $z$ increases from the de-localized productive phase to
the finitely generated phase, and then to the reactive phase. For smaller values of $p$ (typically, 
$|{\cal F}|p\lesssim 10^{-2}$, so that $\sqrt{|{\cal F}|p}$ is small enough), the numerical
agreement with our $p\to 0$ is good, but it is necessary to zoom out a very small 
neighborhood of $z=1$ to observe the transition.  

\begin{figure}[t]
\centering
\includegraphics[scale=0.6]{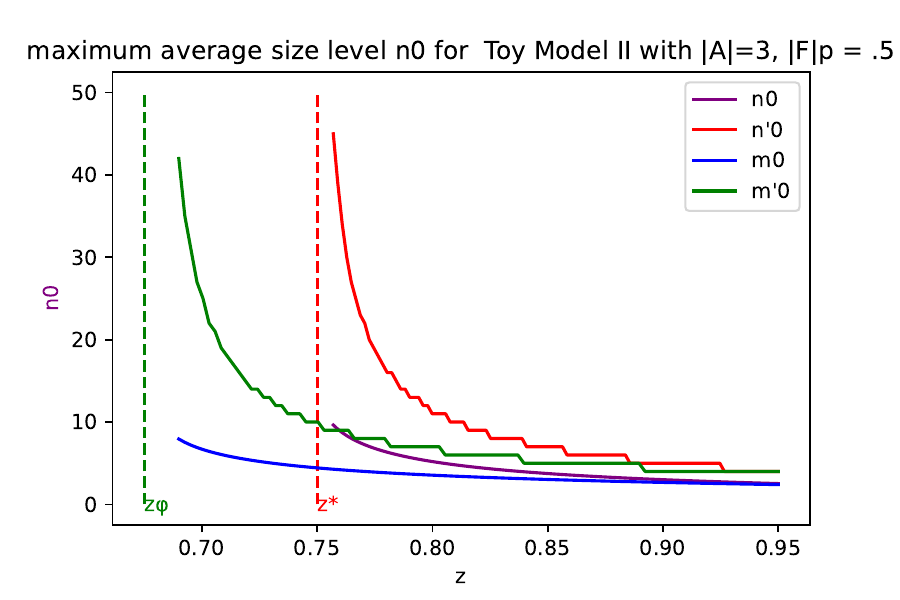}  
\caption{Levels $n_0,n'_0$ of Model II in the reactive phase $z>0.746...$, and primitive  levels $m_0,m'_0$ 
in finitely generated and reactive phases, $z>0.676...$ (parameters : $|{\cal A}|=3$ and  $|{\cal F}|p = 0.5$)}
\label{fig:n0}
\end{figure}

\begin{figure}[t]
\centering
\includegraphics[scale=0.6]{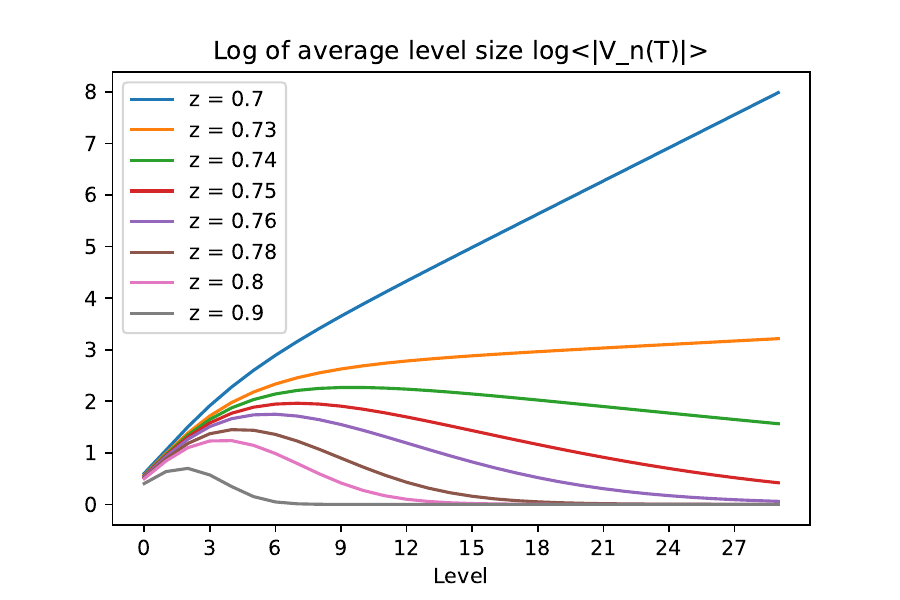}  
\caption{Log of average population $\log \langle |V_n(\T)| \rangle $ in anabolic Model II for 
$|{\cal A}|=3$ and $|{\cal F}|p=0.5$ (theoretical predictions).}
 \label{fig:log-pop-II-z}
\end{figure}

\begin{figure}[t]
\centering
\includegraphics[scale=0.6]{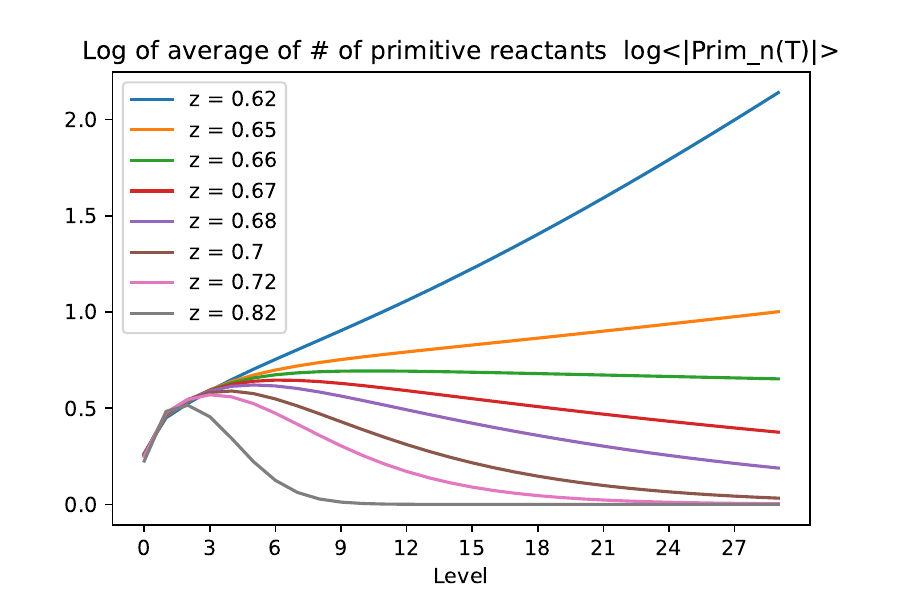}
\caption{Log of average number of primitive  reactants $\log \langle |\Prim_n(\T)| \rangle $ in anabolic Model II for 
$|{\cal A}|=3$  and $|{\cal F}|p=0.5$ (theoretical predictions).}
 \label{fig:log-prim-II-z}
\end{figure}

\Medskip 
\noindent {\em Average population and average number of primitive reactants.}  We plot in Figure 
\ref{fig:log-pop-II-z}, resp. Figure \ref{fig:log-prim-II-z}, the prediction
of our model for $\log \langle |V_n(\T)| \rangle $ using (\ref{eq:II-Vn}),  resp. for $\log \langle |\Prim_n(\T)| \rangle $ using (\ref{eq:II-Primn}). In both cases, we have assumed for definiteness that 
$|F|=1$ for all $F\in {\cal F}$. The phase change at $z=z^*$, resp. $z=z_{\phi}$, is quite apparent.

%%%%%%%%%%%%%%%%%
%%%%%%%%%%%%%%

\section{Conclusion and perspectives}

%%%%%%%%%%%%%%%%%%
%%%%%%%%%%%%%%%%%

We have dealt within the limited scope of this work with a family of "Small-World" chemical reaction models 
(RRBCN Model I, generalized as RRBCN Model II) dependent on an atom set,
a foodset and a fugacity parameter $z$ conjugate to the complexity index. 
 RRBCN networks were explored along two directions:
 
\textbullet\ dynamical patterns involving autocatalysis and multistability were systematically
investigated and found in the low levels of networks (up to levels 5-6, i.e. 6-7 atoms). Due to 
CPU time and memory limitations, it seems difficult
to extend simulations to much longer molecules. 

\textbullet\ statistics of individual mechanisms were established by using an approximate mapping to  an inhomogeneous Galton-Watson model. Our predictions agree surprisingly well with numerics.  Whereas the model predicts only a finite (though potentially large) number of primitive reactions when reactivity is independent of the complexity index  (Model I, $z=1$), a number of qualitatively different phases is observed otherwise (Model II). In particular, in an intermediate region for $z$ far from 0 and 1, there appears a productive phase, characterized by an infinite number of mechanisms with increasing 
complexity.

\Medskip The inspiration for the dynamical patterns looked for in \S \ref{section:low-level} 
comes in part from the above cited article \cite{BunRiv}, where the authors introduced a toy
model of networks
of  inhibitory interactions based on a sparse graph \textsf{G}  playing the same role for dynamics as \textsf{ConnG} (stationary states are indexed by independent subsets of \textsf{G}). Stochastic dynamical simulations displayed very interesting evolutionary features characteristic
of glassy dynamics, see also \cite{Bar} for a study along those lines of the underlying hard-core
model on random graphs. It would thus be interesting to see if large, sparse \textsf{ConnG} graphs appear by raising the cut-off level $n_{max}$, and if conclusions of the above article 
are also relevant for RRBCN models. 
   
\Bigskip Despite these achievements, our model is a toy model of reactivity, because it fails to deal with the following points:

\begin{enumerate}
\item a larger number of kinetic rankings, depending mainly on the order of magnitude of the kinetic
rate, should be distinguished;
\item molecular graphs in biochemistry are seldom linear. Instead, they branch, form (e.g. aromatic) cycles, 
and include double and triple bonds, respecting valence rules; 
\item reactions involve electronic rearrangements and protonation/deprotonation steps which
allow bond formation/breaking. Thus the product of a reaction is not obtained by pure concatenation/fragmentation of words;    
\item the probability of a bond formation/breaking is strongly dependent on the nature of the atoms involved;
\item anabolic reactions may involve two  reactants, neither of which is 
in the foodset. Furthermore, some more complicated mechanisms (e.g. 2-2 reactions) may be considered.  
\end{enumerate} 

\Medskip However, we have been working on generalizations of our 
model which go a long way towards overcoming these limitations. The principles are the following:

\begin{enumerate}
\item we need a reasonable  prior for kinetic parameters, based on chemical expertise and computational chemistry \cite{NanUntStuNgh}; 
\item looking at functional molecules relevant in extant biology (see e.g. \cite{McMBeg})  and prebiotic experiments 
 makes it clear that an overwhelming proportion of molecules is essentially in the form of a linear skeleton, with some short and specific cycles  (which can be thought of as additional atoms in our language), and branches of length one (side methyl groups,...), provided functional groups such as carbonyl, carboxyl, etc. are considered as atoms.  Thus the 1-derivation model
 (see Suppl. Info., {\tt \S III.B.2}) captures most relevant
molecular structures;
\item rearrangements can by-and-large be predicted by the  "arrow-pushing diagrams" of undergraduate biochemistry, which change our rules only marginally;
\item the dependence of kinetic parameters on the atoms makes part of the  prior in 1., and 
can be handled by considering a mapping to a multi-type Galton-Watson model, and using large-deviation theoretic tools to deal with the energy landscape.  
\end{enumerate}

Thus only point 5. remains, for which we have no particular strategy. Even before trying to cope with it, much remains to be said about
the structure of the reaction network in our toy model. Here only a few preliminary results about
the anabolic reaction network, and a remark about the catabolic reaction network, were added in Suppl. Info.,
  showing the relevance of the auxiliary tree model. We already mentioned the difficulty
  to superpose the two parts of the network (work in progress).  
But many interesting formulas (not included) can be obtained in various asymptotic regimes. This raises 
very intriguing questions.  There is a small number of chemically interpretable  "hyperparameters" in the model, such as $p,q$; 
$z$ (or a level-dependent version of $z$); physico-chemical parameters (solvent, pH, temperature...) which change the energy landscape, etc.  The first question is: can one 
infer the value of hyperparameters from experimental data ? Which phase is relevant ? Can one 
discuss the plausibility of a phase transition as one of the hyperparameters is varied ?  Does one
observe different statistical behaviors of reaction networks depending on the level, which would
point out at a possible chemical evolution towards subnetworks involving longer molecules ?

\Bigskip 
\subsection*{Funding}  
The author did not receive support from any organization for the submitted work.

\subsection*{Conflicts of interest}  
The author has no relevant financial or non-financial interests to disclose.

\subsection*{Acknowledgement}
We thank Olivier Rivoire for insightful discussions and helpful comments on the manuscript.

\subsection*{Code availability}  
Scripts used to reproduce the analyses and generate the figures reported in this article can be provided upon request to the corresponding author.

%%%%%%%%%%%%%%%%%%%%%%%%%%%%%%%%%%%%%%%%%%%%%%%%%%%


\begin{thebibliography}{100}

\bibitem{Adam} Z. R. Adam, A. C. Fahrenbach, S. M. Jacobson, B. Kacar, D. Y. Zubarev (2021). {\em Radiolysis generates a complex organosynthetic chemical network}, Sci. Rep. {\bf 11}, 1743.

\bibitem{AlbBar} R. Albert, A.-L. Barabasi  (2002). {\em Statistical mechanics of complex networks}, Rev. Mod. Phys. {\bf 74 (47)}, 47–97.


\bibitem{And0} 
J. L. Andersen, C. Flamm, D. Merkle, and P. F. Stadler (2016). {\em  A software
package for chemically inspired graph transformation}, in:  Graph
Transformation, ICGT 2016, ser. Lecture Notes Comp. Sci. {\bf 9761}, R. Echahed
and M. Minas, Eds., Springer.

\bibitem{And1}  J. L. Andersen, C. Flamm, D. Merkle, and P. F. Stadler (2017).  {\em An intermediate level of
abstraction for computational
systems chemistry},  Phil. Trans. R. Soc. {\bf A 375}, 20160354. 




\bibitem{And2}  J. L. Andersen, C. Flamm, D. Merkle, P. Stadler (2019). {\em Chemical Transformation Motifs -- Modelling Pathways as Integer Hyperflows},  IEEE/ACM Transactions on Computational Biology and Bioinformatics {\bf 16 (2)}. 

\bibitem{AthNey}  K. B. Athreya, P. E. Ney (1972). {\em Branching processes}, Grundlehren der 
Mathematischen Wissenschaften {\bf 196},  Springer.

\bibitem{Bar} J.Barbier, F. Krzakala et al. (2013). {\em  The hard-core model on random graphs revisited}, J. Phys. Conf. Ser. {\bf 473}, 012021.

\bibitem{BloLacNgh} A. Blokhuis, D. Lacoste, P. Nghe (2020). {\em Universal motifs and the diversity of autocatalytic systems}, PNAS  {\bf 117}, 25230-25236.


\bibitem{BraSmi} R. Braakman, E. Smith (2013). {\em The compositional and evolutionary logic of metabolism },  Phys. Biol. {\bf 10},  011001.

\bibitem{Bre} R. Breslow. {\em On the mechanism of the formose reaction},
Tetrahedron Lett. {\bf 21}, 22–26 (1959).

\bibitem{BunRiv} G. Bunin, O. Rivoire (2025). {\em Evolutionary features in a minimal physical system: Diversity, selection, growth, inheritance, and adaptation}, PNAS {\bf 122} (31). 

\bibitem{But} Butlerow, A. {\em Formation synthétique d’une substance sucrée}, C. R.
Acad. Sci. 53, 145–147 (1861).


\bibitem{ChemKin} CHEMKIN, chemistry simulation software, https://www.ansys.com/fr-fr/products/fluids/ansys-chemkin

\bibitem{CraFei} G. Craciun, M. Feinberg (2016). {\em  Multiple Equilibria in Complex Chemical
Reaction Networks: II. The Species-Reaction Graph}, SIAM J Appl Math.
{\bf 66}, 1321–1338. 

\bibitem{EcoCyc} EcoCyc E. Coli database, https://ecocyc.org/.

\bibitem{Eigen} M. Eigen (1971). {\em Selforganization of Matter
and the Evolution of Biological Macromolecules}, Naturwissenschaften {\bf 58 (10)}.  

\bibitem{ErdRen}  P. Erd\"os, A. R\'enyi  (1960). {\em On the evolution of random graphs}, Publ.  Math. Inst.  Hung. Acad.  Sci.  {\bf 5}, 17–61. 

\bibitem{FlaHel} C. Flamm, M. Hellmuth, D. Merkle et al. (2020).  {\em Generic Context-Aware Group Contributions},  IEEE/ACM Transactions on
Computational Biology and Bioinformatics.
 
\bibitem{GagBlaSmiBau} P. Gagrani, V. Blanco, E. Smith, D. Baum (2024). {\em Polyhedral geometry and combinatorics of an autocatalytic ecosystem}, J. Math. Chem. {\bf 6}2, 1012–1078.

\bibitem{Gho} G. Ghoshal, V. Zlatic, G.Caldarelli, M. E. J. Newman (2009). 
{\em  Random hypergraphs and their applications},  Phys Rev {\bf E79}, 066118. 




\bibitem{Hor} W. Hordijk (2019). {\em A History of Autocatalytic Sets.
A Tribute to Stuart Kauffman}, Biological Theory {\bf 14}, 224–246. 


\bibitem{HorSte} W. Hordijk, M. Steel (2004). {\em Detecting autocatalytic, self-sustaining sets in chemical reaction systems},  J. Theor. Biol. {\bf 227}, 451–461.

\bibitem{IsmChaHab}
Idil Ismail, Raphael Chantreau Majerus, and Scott Habershon
The Journal of Physical Chemistry A 2022 126 (40), 7051-7069
DOI: 10.1021/acs.jpca.2c06408

\bibitem{Joshi} B. Joshi, A. Shiu (2013). {\em Atoms of multistationarity in chemical reaction networks}, J. Math. Chem. {\bf 51}, 153-178. 

\bibitem{Kauf} S. A. Kauffman (1971). {\em Cellular Homeostasis, Epigenesis and Replication in Randomly Aggregated Macromolecular Systems.},  J. Cybern. {\bf 1}, 71–96. 


\bibitem{KEGG} KEGG (Kyoto Encyclopedia of Genes and Genomes),  database resource for understanding high-level functions and utilities of the biological system, see https://www.genome.jp/kegg/.   


\bibitem{KobSug} T. J. Kobayashi and Y. Sughiyama (2017). {\em Stochastic and information-thermodynamic structures of population dynamics
in a fluctuating environment}, Phys. Rev. E {\bf 96}, 012402.

\bibitem{Mar}
Margraf, J.T., Jung, H., Scheurer, C. et al. Exploring catalytic reaction networks with machine learning. Nat Catal 6, 112–121 (2023). https://doi.org/10.1038/s41929-022-00896-y

\bibitem{MOD} M\O D, a software package developed for graph-based cheminformatics, https://cheminf.imada.sdu.dk/mod/.


\bibitem{MulFlaSta} S. Müller, C. Flamm, P. Stadler (2002). {\em What makes a reaction network “chemical”?}, Journal of Cheminformatics, 14:63.









\bibitem{McMBeg} J. Mc Murry, T. P. Begley (2005). {\em The organic chemistry of biological pathways}, Roberts and Company Publishers. 







\bibitem{Pan} E. de Panafieu  (2015). {\em  Phase transition of random non-uniform hypergraphs},
 J. Discrete Alg. {\bf 31}, 26–39. 

\bibitem{Peng} Z. Peng Z, J. Linderoth, D. A. Baum  (2022). {\em The hierarchical organization of autocatalytic reaction networks and its relevance to the origin of life}, PLoS Comput. Biol. 
{\bf 18(9)}: e1010498. 

\bibitem{PubChem} PubChem chemical database, https://pubchem.ncbi.nlm.nih.gov/.



\bibitem{Rav} E. Ravasz, A. L. Somera, D. A. Mongru, Z. N. Oltvai, A.-L. Barabasi (2002).   {\em Hierarchical Organization of
Modularity in Metabolic
Networks}, Science {\bf 297}. 

\bibitem{RMG} RMG, Reaction Mechanism Generator, https://rmg.mit.edu/.

\bibitem{Huck} W. E. Robinson, E. Daines, P. van Duppen, T. de Jong, W. T. S. Huck  (2022). {\em Environmental conditions drive self-organization of reaction pathways in a prebiotic reaction network}, Nature Chemistry {\bf 14}, 623–631.

\bibitem{Mar}
Margraf, J.T., Jung, H., Scheurer, C. et al. (2023). Exploring catalytic reaction networks with machine learning. Nat Catal 6, 112–121. https://doi.org/10.1038/s41929-022-00896-y.


\bibitem{NanUntStuNgh} P. Nandan, J. Unterberger, T. Stuyver, P. Nghe. {\em Universality of carbon metabolism in Reaction
networks of small molecules}, in preparation.

\bibitem{UnsGriRei} Jan P. Unsleber, Stephanie A. Grimmel, and Markus Reiher
Journal of Chemical Theory and Computation 2022 18 (9), 5393-5409
DOI: 10.1021/acs.jctc.2c00193.

\bibitem{RufDan}  A. Ruf, G. Danger (2022). {\em Network Analysis Reveals Spatial Clustering and Annotation of
Complex Chemical Spaces: Application to Astrochemistry}, Anal. Chem. {\bf 94}, 14135-14142.

\bibitem{SteHor} M. Steel, W. Hordijk, J. C. Xavier (2019). {\em Autocatalytic networks in biology: structural theory and algorithms}, J. R.
Soc. Interface {\bf 16}, 20180808.

\bibitem{UntNgh} 	P. Nghe, J. Unterberger (2022). {\em Stoechiometric and dynamical autocatalysis for diluted chemical reaction networks}, J. Math. Biol. {\bf 85 (26)}.


\bibitem{WagFel} A. Wagner, D. A. Fell (2001). {\em The small world inside large metabolic networks}, 
Proc. R. Soc. Lond. {\bf B268}, 1803-1810. 

\bibitem{Wang}
Wang, LP., Titov, A., McGibbon, R. et al. (2014). Discovering chemistry with an ab initio nanoreactor. Nature Chem 6, 1044–1048. https://doi.org/10.1038/nchem.2099.


\bibitem{WatStr} D. J. Watts,  S. H. Strogatz (1998). {\em  Collective dynamics of
"small-world" networks}, Nature {\bf 393}.

\bibitem{Wei} D. Weininger (1988). {\em SMILES, a chemical language and information system. 1. Introduction
to methodology and encoding rules}, J.  Chem. Info.  Comp. Sci. {\bf  28:1} (1988), 31-36.

\bibitem{WenSpo} Wen, M., Spotte-Smith, E.W.C., Blau, S.M. et al. Chemical reaction networks and opportunities for machine learning. Nat Comput Sci 3, 12–24 (2023). https://doi.org/10.1038/s43588-022-00369-z

\bibitem{ZhaGarSav}
Qiyuan Zhao, Sanjay S. Garimella, and Brett M. Savoie (2023).
Journal of the American Chemical Society 145 (11), 6135-6143.


\bibitem{ZhaSav}
Zhao, Q.; Savoie, B. M. Simultaneously improving reaction coverage and computational cost in automated reaction prediction tasks. Nat. Comput. Sci. 2021, 1, 479– 490,  DOI: 10.1038/s43588-021-00101-3

\end{thebibliography}
\end{document}